\documentclass[11pt]{article}

\usepackage{graphicx,amsmath,upgreek,bm}
\usepackage{colortbl}
\usepackage[rgb, table, dvipsnames]{xcolor}
\usepackage{booktabs}
\usepackage{makecell}
\usepackage{float}
\usepackage{pifont}
\usepackage{todonotes}
\usepackage[toc,page]{appendix}
\usepackage{multirow}
\usepackage{subcaption} 	
\usepackage{amsmath,amsfonts,amsthm,bm,mathtools}
\usepackage{lipsum}
\usepackage{geometry}
\usepackage{algorithm}
\usepackage[noend]{algpseudocode}
\usepackage{custom_format}  %
\usepackage{arydshln}
\usepackage{quotes}
\makeatletter

\renewcommand*\env@matrix[1][*\c@MaxMatrixCols c]{%
\hskip -\arraycolsep
\let\@ifnextchar\new@ifnextchar
\array{#1}}
\makeatother
\definecolor{refcolor}{rgb}{0.65,0,0.35}
\usepackage[colorlinks=black]{hyperref}
\hypersetup{
	colorlinks = true,
	linkcolor = refcolor,
	citecolor = refcolor,
	filecolor = refcolor,      
	urlcolor = blue,
	citecolor = refcolor
}

\graphicspath{ {figures/}{../figures/} }
\begin{document}
\title{\textbf{A Topological Framework for Identifying Phenomenological Bifurcations in Stochastic Dynamical Systems}}
\author{Sunia Tanweer \and Firas A.~Khasawneh \and Elizabeth Munch \and Joshua R. Tempelman}
\date{Dept. of Mech. Engineering, Michigan State University}

\maketitle

\begin{abstract}

Changes in the parameters of dynamical systems can cause the state of the system to shift between different qualitative regimes. These shifts, known as bifurcations, are critical to study as they can indicate when the system is about to undergo harmful changes in its behavior. In stochastic dynamical systems, there is particular interest in P-type (phenomenological) bifurcations, which can include transitions from a mono-stable state to multi-stable states, the appearance of stochastic limit cycles, and other features in the probability density function (PDF) of the system's state. Current practices are limited to systems with small state spaces, cannot detect all possible behaviours of the PDFs, and mandate human intervention for visually identifying the change in the PDF. In contrast, this study presents a new approach based on Topological Data Analysis (TDA) that uses superlevel persistence to mathematically quantify P-type bifurcations in stochastic systems through a "homological bifurcation plot''---which shows the changing ranks of $0th$ and $1st$ homology groups. Using these plots, we demonstrate the successful detection of P-bifurcations on the stochastic Duffing, Raleigh-Vander Pol and Quintic Oscillators given their analytical PDFs, and elaborate on how to generate an estimated homological bifurcation plot given a kernel density estimate (KDE) of these systems by employing a tool for finding topological consistency between PDFs and KDEs.

\end{abstract}

\section{Introduction}

Determination of the state of a stochastic dynamical system and the inception of a bifurcation is of immense interest to accurately understand, predict and control \cite{Annunziato2010, ANNUNZIATO2013487, Mifeng2019} the behaviour of the system. In such systems, bifurcations may be induced by a change in system properties such as the damping or stiffness \cite{Mamis2016, ZHAO2016166, LIU2020103085} or excitation properties such as the intensity, color, and bandwidth of forcing noise \cite{Arnold1995,Zakharova2010,Falk2018,Zakharova2010}. Additionally, phenomenon such as noise bursts have also been shown to induce bifurcations in the system \cite{Falk2017}. Since stochastic dynamics have seen wide use across an array of applications including chemical reactions \cite{Falk2017,Schlogl1972}, fluid dynamics \cite{VENTURI2010}, vibration control \cite{Liu2018}, energy harvesting \cite{McInnes2008, tai2020optimum}, genetics \cite{Lee2018}, and financial markets \cite{Chiarella2008}, there have been many thrusts to qualitatively and quantitatively analyze stochastic bifurcations from a variety of academic communities \cite{Falk2018,arnold2013dynamical,Liu2018,Mamis2016,Schenk-Hoppe1996a}.

In stochastically excited systems, bifurcations are characterized in two distinct types: P (Phenomenological) and D (Dynamical). While the D-type can be quantitatively measured using sudden changes in the sign of the largest Lyapunov exponent, the P-type bifurcation is a qualitative change in the topology of the joint probability density function (referred to as PDF onwards) of the system response. Based on this, a P-bifurcation would be, for example, the transition between mono-stable and bi-stable PDFs \cite{Yang2016,Yang2017,Schlogl1972,Falk2017,ARNOLD1996,Arnold1995}, whereas a D-type bifurcation would be a separation of state space regions between which no transitions may occur \cite{Falk2018,Arnold1995, Schenk-Hoppe1996a,ARNOLD1996,li2019}. It is noteworthy that the P- and D-types are unrelated phenomenon \cite{xu_gu_zhang_xu_duan_2011} and may not occur in the same regimes \cite{kumar_narayanan_gupta_2016}. 

Given that P-bifurcation are attributed to the qualitative changes in the PDF, their analysis has traditionally been restricted to qualitative methods \cite{KUMAR2022104086}. In \cite{Falk2017}, fixed points in the PDF are computed for a one dimensional stochastic bi-stable system. These maxima and minima correspond to favorable and unfavorable points, respectively, which are used to show where in the phase portrait dynamics tends to. In \cite{Falk2018}, the fixed points of the system PDF are used with information from limit cycles in the stochastic phase portrait to indicate ridges in the PDF. In \cite{Song2010}, a stochastic bi-stable system's bifurcations are captured in an empirical bifurcation diagram similar to those presented in \cite{Ozbudak2004}. The resulting two dimensional empirical bifurcation diagram uses density estimators to capture the trajectory distribution at different parameter values, and then these distributions are plotted on a vertical axis which varies on the parameter space. This results in a qualitative discernment of mono-stable versus bi-stable behavior. Xu et al. \cite{Xu2011} also use extrema of the PDF from an exact stationary solution to the stochastic Duffing-Van der Pol equation to show bistable vs monostable behaviour. However, unlike \cite{Song2010,Ozbudak2004}, the analytic evaluation of \cite{Xu2011} allows for the exact determination of the number of modes in the stationary distribution. Thus, \cite{Xu2011} constructs a two-dimensional bifurcation diagram by varying parameters, tracking the number of modes in the exact distribution, and the resulting bifurcation diagram displays system parameters on the horizontal and vertical axes and shades in regions according to the number of modes found at that parameter pair. This procedure is utilized to construct the bifurcation diagrams in \cite{Li2017,li2019coherence}. 

While the methods provided in \cite{Falk2017,Falk2018,Xu2011,Song2010,Ozbudak2004} offer decent qualitative approaches for characterizing stochastic bifurcations, there are limiting factors to the current state-of-the-art. For instance, these measures only work for one dimensional systems. If a joint-distribution is assumed, it would be much more challenging to display the shifting dynamics as is done by \cite{Song2010,Ozbudak2004} in higher dimension. Additionally, simply picking the number of peak values in the PDF would fail to capture richer dynamical behavior such as limit-cycles \cite{Mamis2016}. Studies which have approached bifurcations of multi-dimensional systems typically display multiple joint distributions across a span of parameters and draw qualitative conclusions to claim within which range of parameters a P-bifurcation occurs in \cite{Mamis2016,Djeundam2013, Liu2018}. Lastly, a major challenge associated with these qualitative methods is their heavy reliance upon visual inspection of the PDF, making it difficult to estimate the stability bounds for the system \cite{kumar_narayanan_gupta_2016, kumar_narayanan_gupta2_2016}. Some recent attempts have been made to quantify the onset or occurrence of a P-type bifurcation using changes in the sign of Shannon entropy, but the technique presented is not applicable to all generic nonlinear dynamical systems and bifurcation detection may be difficult in systems with multi-dimensional state space \cite{venkatramani_sarkar_gupta_2018, kumar_narayanan_gupta2_2016}. 

Recognizing this dearth of quantitative techniques for identifying the topology and subsequently the changes in topological features of the PDF, our primary focus in this paper is proposing a technique for quantification and characterization of P-bifurcations. We provide a quantitative tool for analyzing the qualitative P-bifurcation which works for both the analytical probability distributions generated by the limited number of solvable SDEs, and also extends to systems with no known exact solution. Specifically, we apply persistent homology, a popular tool from Topological Data Analysis \cite{TamalDey2022}, to PDFs or Kernel Density Estimates (KDE) of the system response. This allows us to quantify the shape of the probability distributions which provides a quantitative bifurcation diagram where the change in the ranks of various homology classes in the topology of the distribution identify parameter values which induce a shift in the system dynamics. 

Consider the general form used to describe an n-dimensional stochastic dynamical system:
\begin{equation}
d\textbf{X} = \boldsymbol{\mu}(\textbf{X}, t) dt + \boldsymbol{\sigma}(\textbf{X}, t)d\textbf{W}(t)
\end{equation}
where d\textbf{W}(t) is an n variable Wiener process, and $\boldsymbol{\mu}$ and $\boldsymbol{\sigma}$ are drift and diffusion function vectors respectively, such that $\boldsymbol{\sigma}$ induces parametric stochastic excitation. In this work, the analytical PDFs used for systems have not been derived directly. It has been assumed that either the analytical PDF is available, is easily derivable or if not possible, a KDE can be numerically computed using system's response. However, the eager reader can refer to methods of equivalent non-linear system \cite{xuwei_yang_2016, CAI1988409, CAUGHEY19862}, stochastic averaging \cite{Roberts1986,Zhu1996,Xu2011}, perturbation for determining stochastical moments \cite{KAMINSKI2010272, bonizzoni_nobile_2014, williams_larsen_2001} and differential transform method \cite{Hesam2012} for more information. Such stochastic differential equations which take the Ito form may also be split into the Fokker-Plank (FP) equation, which is an approximation of the Kramers-Moyal expansion and brings the SDE to the form \cite{gardiner_2002}
\begin{equation}
\partial_{t} p(\textbf{X}, t) = -\sum_i \partial_{i} [\mu_{i}(\textbf{X}, t)p(\textbf{X}, t)] + \frac{1}{2}\sum_{i,j} \partial_i \partial_j [\{\boldsymbol{\sigma}(\textbf{X}, t)\boldsymbol{\sigma}^\textbf{T}(\textbf{X}, t)\}_{ij}p(\textbf{X},t)]
\end{equation}
where $p(\textbf{X},t)$ is the joint probability distribution of the state variables. In \cite{Mamis2016,Caughey1982, Bezen1996,Yong1987}, exact stationary solutions are achieved for the reduced Fokker-Plank-Kolormogrov (rFPK) equation (same as FP but with a stationary joint distribution of the state variables). Using the mentioned tools, past research has qualitatively tracked the bifurcation behavior of a small subset of stochastic differential equations. Nevertheless, it is essential to note that the the source of the probability distribution is inconsequential to the methods proposed in this paper.

\section{Mathematical Preliminaries}
\label{sec:maths}

This work relies on a range of mathematical tools from different fields, including topological data analysis (TDA) \cite{TamalDey2022}.
The key mathematical technique used here is persistent homology computed using both simplicial and cubical complexes. 
In addition, we also require the use of KDEs and employ a tool for determining their topologically consistent persistent homology. 
This section presents a brief introduction to each of these techniques.

\subsection{Topological Spaces}

The first topological space we will use is that of a cubical complex, which is useful in the case of discretized function (matrix) data. 
A cubical complex is built by elementary intervals, defined as either a unit interval $[u, u+1]$, or a degenerate interval $[u,u]=[u]$ where $u\in \mathbb{Z}$. 
Then, for an $n$-dimensional space, a cube is defined by the product $I = \prod\limits_{i = 1}^{K}[u_i, u'_i]$ where $u'_i \in \{u_i, u_{i+1}\}$. The dimension of the resulting cube is the number of non-degenerate intervals in this product. For example, $0$-cubes are vertices, $1$-cubes are edges, $2$-cubes are squares, while $3$-cubes are voxels, and so on. 
If one cube $\sigma$ is a subset of another $\tau$, we say $\sigma$ is a face of $\tau$ and write $\sigma \leq \tau$. 
For example, the edges $\sigma_1 = [0] \times [0,1]$ and $\sigma_2 = [0,1] \times [0]$ are both faces of the square $\tau = [0,1]\times [0,1]$. 
A cubical complex, then, is a collection of cubes where including a $d$-cube in the complex implies that all of its faces are included. 

Another form of complex we utilize in this work is that of a simplicial complex. 
Given a collection of vertices $V$, a $d$-simplex is a set of $d+1$ vertices, $\sigma = [v_0,\cdots, v_d]$. 
In this case, a 0-simplex is still a point and a 1-simplex is still an edge, however, a 2-simplex is a triangle, a 3-simplex a tetrahedron, etc. 
Again, we say that $\sigma$ is a face of $\tau$ if $\sigma \subseteq \tau$, and write $\sigma \leq \tau$. 
A simplicial complex $K$ is a collection of simplices such that if $\sigma \in K$ and $\tau \leq \sigma$, then $\sigma \in K$. 

In the next sections, we give a brief introduction to homology, super- and sublevel-set persistent homology, and direct the interested reader to \cite{TamalDey2022,Steve2015} for full details. 

\subsection{Homology}

For a fixed topological space $X$, which for the purposes of this paper will be a finite combinatorial structure such as a simplicial or cubical complex, we can define homology as follows. 
The $p$-th chain group $C_p(X)$ is a vector space given by finite formal sums of the $p$-dimensional simplices (or cubes), $\sum \alpha_i\sigma_i$, where addition is given by addition of the coefficients. 
For this work, we assume the coefficients $\alpha_i$ are taken from $\Z_2 = \{0,1\}$.
The boundary map takes a simplex to the sum of its faces with dimension exactly one lower: $\partial_p: C_p(x) \xrightarrow{} C_{p-1}(x)$ %
Then homology is defined as $H_p(X) = \Ker(\partial_p) / \im(\partial_{p+1})$.
In particular, an element of the $p$-dimensional homology is represented by a collection of $p$-simplices such that the boundary map applied to it is 0, and where if it is equivalent to another representative if they differ by the boundary of a higher dimensional collection of simplices.
More intuitively, the type of structure being measured depends on the dimension of homology being used --- connected components are captured by zero-dimensional homology, loops by one-dimensional homology, and voids by two-dimensional homology. 

\subsection{Persistent Homology}
Persistent homology is a widely used technique in TDA that captures information about the shape of a parameterized space by tracking how its homology changes as the parameter varies \cite{TamalDey2022}. 

Say we are given a topological space with a function $f:X \to \R$. 
For a collection of real values, $a_1\leq a_2 \leq \cdots \leq a_n$, we have a \textit{sub-level-set filtration} given by 
\begin{equation*}
    X_{a_1} \subseteq X_{a_2} \subseteq \cdots \subseteq X_{a_n}
\end{equation*}
where $X_a = f\inv(-\infty,a]$. 
Similarly, there is a \textit{super-level-set filtration} 
\begin{equation*}
    X^{a_n} \subseteq X^{a_{n-1}} \subseteq \cdots \subseteq X^{a_1}
\end{equation*}
where $X^a = f\inv[a,\infty)$. 
It is a fundamental theorem in algebraic topology that homology is functorial; that is to say that if we have a map from one space to another (in our case inclusions $X \subset Y$), there is an induced map (in this case, a linear map on vector spaces) on the relevant homology groups $H_p(X) \to H_p(Y)$. 
We use this to create a sublevelset homology filtration 
\begin{equation*}
    H_p(X_{a_1}) \to H_p(X_{a_2}) \to \cdots \to H_p(X_{a_n})
\end{equation*}
and a superlevelset homology filtration
\begin{equation*}
    H_p(X^{a_n}) \to H_p(X^{a_{n-1}}) \to \cdots \to H_p(X^{a_1}).
\end{equation*}
Another fundamental theorem of persistent homology is that any filtration of this form can be uniquely decomposed into a collection of pairs $(b,d)$. 
In this case, $b=a_i$ indicates the arrival of a new homology class; i.e.~an element of $H_p(X_{a_i})$ (resp. $H_p(X^{a_i})$) not in the image of $H_p(X_{a_{i-1}})$ (resp.~$H_p(X^{a_{i+1}})$). 
Likewise, $d=a_j$ means that class merges with the image of an older class entering $H_p(X_{a_j})$ (resp.~$H_p(X^{a_{j}})$). 

This information is often drawn in a persistence diagram, with the collection of pairs drawn as points in the plane. 
Note that in the case of sublevelset persistence $b<d$ and in the case of superlevelset persistence, $b>d$. 
Thus sublevelset points are always above the diagonal $\Delta = \{ (x,x) \mid x \in \R\}$, while for superlevelset persistence, the points are always below. 
In this paper, we will be particularly focused on persistence in the two cases below. 

\subsubsection{Superlevelset Persistence of Cubical Complexes}
\begin{figure}[!htbp]
    \centering
    \includegraphics[height = 2in]{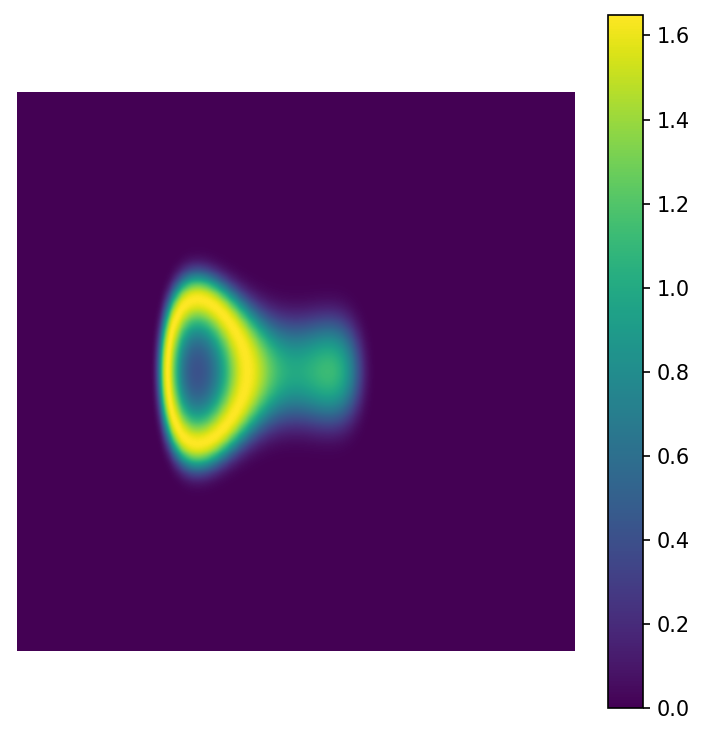} \includegraphics[height = 2in]{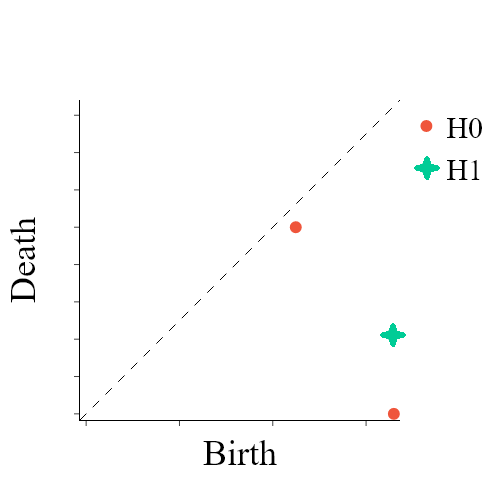}\\
    \includegraphics[height = 1in]{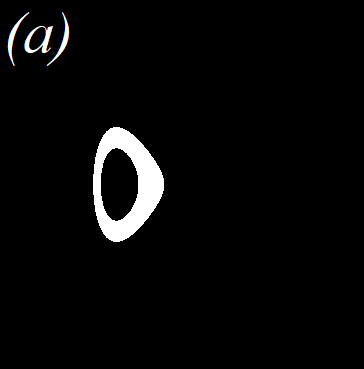}
    \includegraphics[height = 1in]{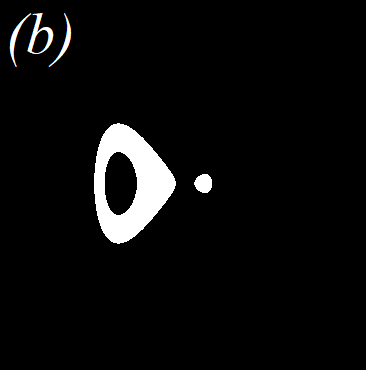}
    \includegraphics[height = 1in]{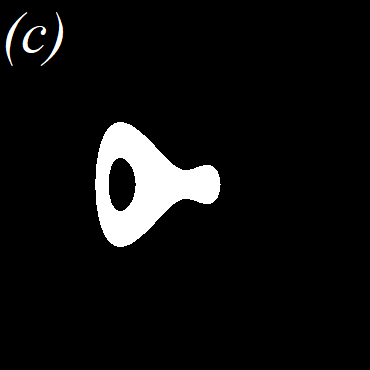}
    \includegraphics[height = 1in]{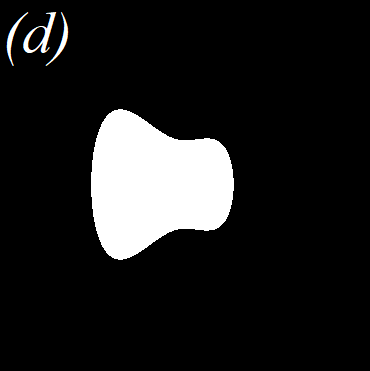}
    \caption{Superlevel cubical persistence of image data. Figures correspond to (a) $K^{1.4}$, (b) $K^{1.05}$, (c) $K^{0.8}$ and (d) $K^{0.2}$. }
    \label{fig:CubicalPersistence}
\end{figure}

Given a grayscale image, or more generally an $m\times n$ matrix $M$, define its domain as the cubical set $K=\mathcal{K}([0,m]\times[0,n])$. %
We can represent the image as a function $f:K \to \mathbb{R}$ defined on this cubical set where for a given cube $s_{i,j}=[i,i+1]\times[j,j+1]$ in $K$ we have $f(s_{i,j})=M_{i,j}$, i.e., the function value is equal to corresponding matrix entry. 
For all lower dimensional cubes $P$ we set $f(P)=\max_{s,j>P} M(s_{i,j})$. 
This way, the superlevelset of $f$ for any function value, which we will denote in this setting as $K^a := f\inv (-\infty,a]$ is also a cubical complex. 
Then we can compute superlevelset persistence of this function. 
See the example in Fig.~\ref{fig:CubicalPersistence}. As $a$ is varied, we get different cubical complexes, $K^{a}$. The connected components and loops in (a-d) are (1, 1), (2, 1), (1, 1) and (1,0). 
\subsubsection{Sublevelset Persistence of Simplicial Complexes from Point Clouds}
\begin{figure}[!htbp]
    \centering
    {\includegraphics[height = 2in]{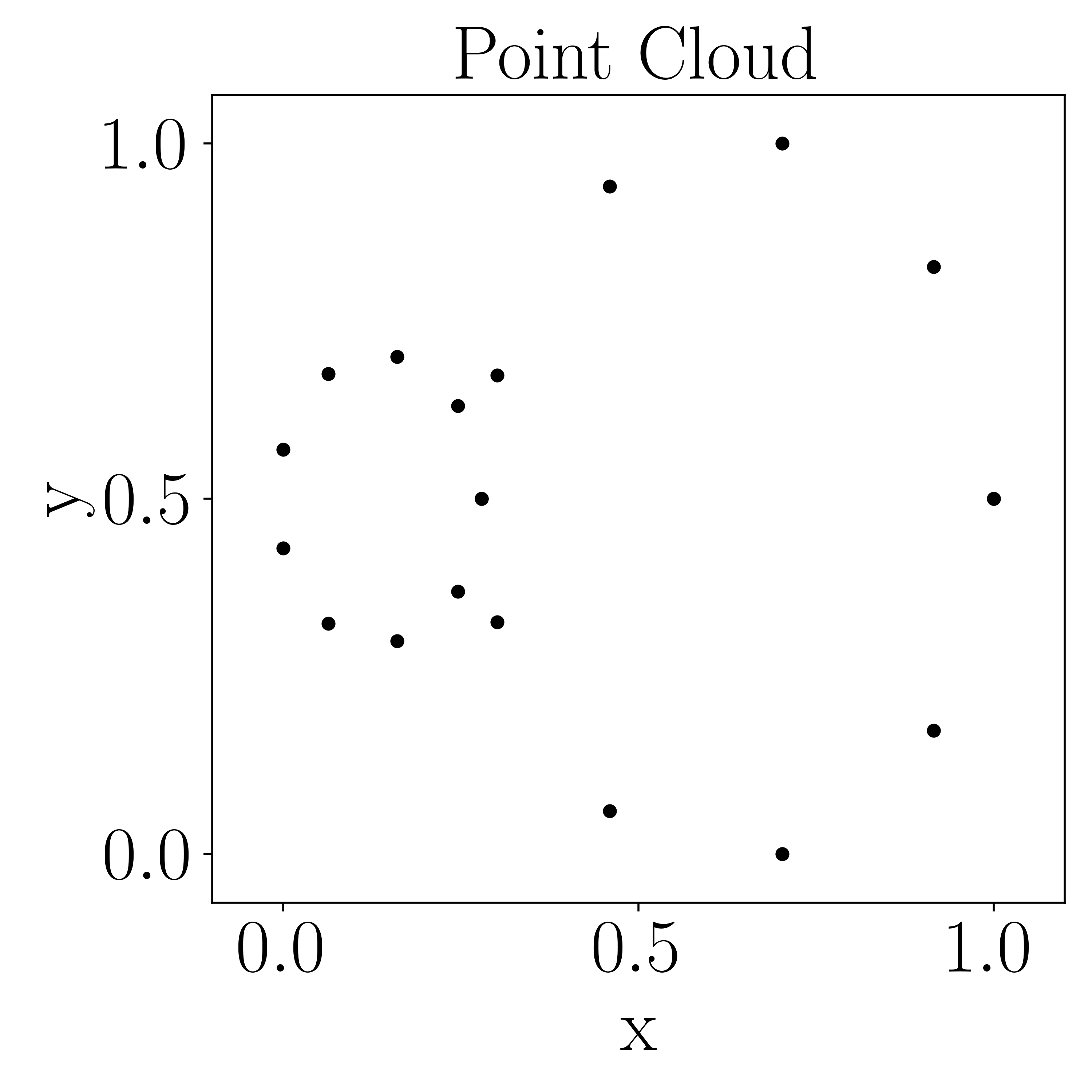}} \includegraphics[height = 2in]{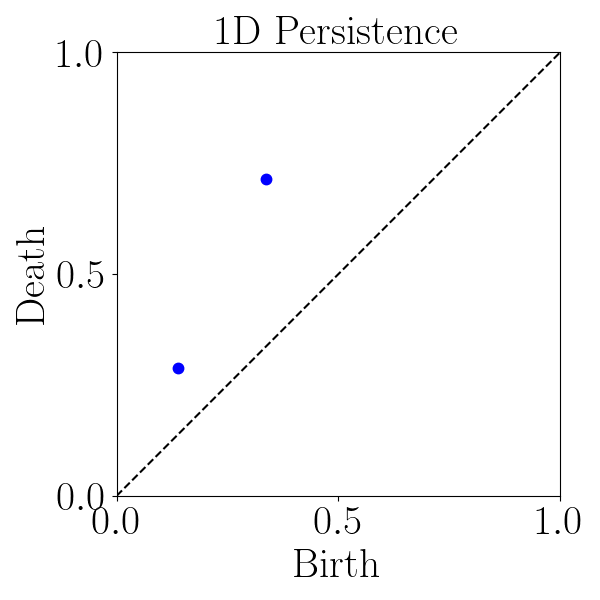}\\
    \includegraphics[height = 1in]{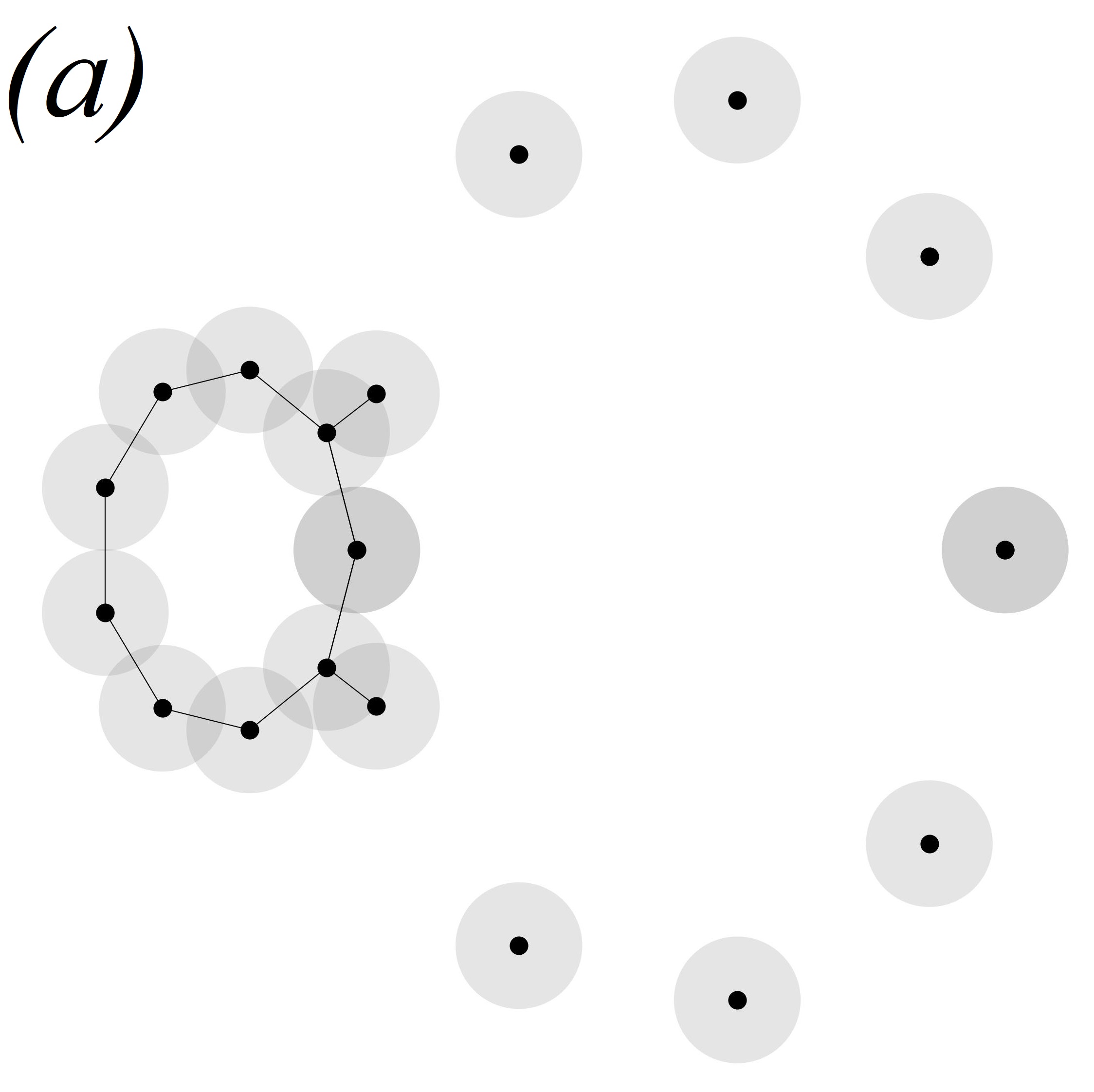}
    \includegraphics[height = 1in]{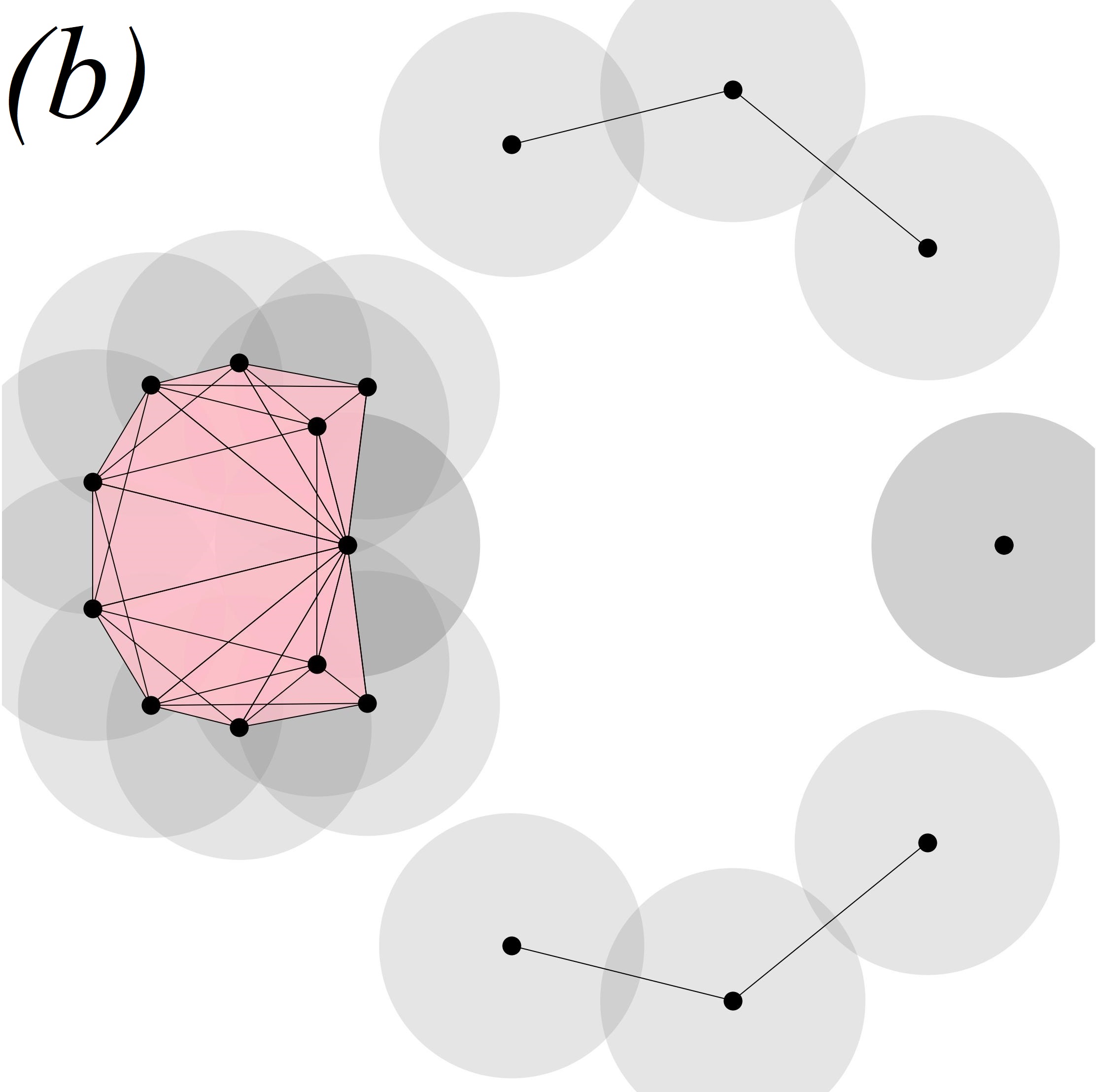}
    \includegraphics[height = 1in]{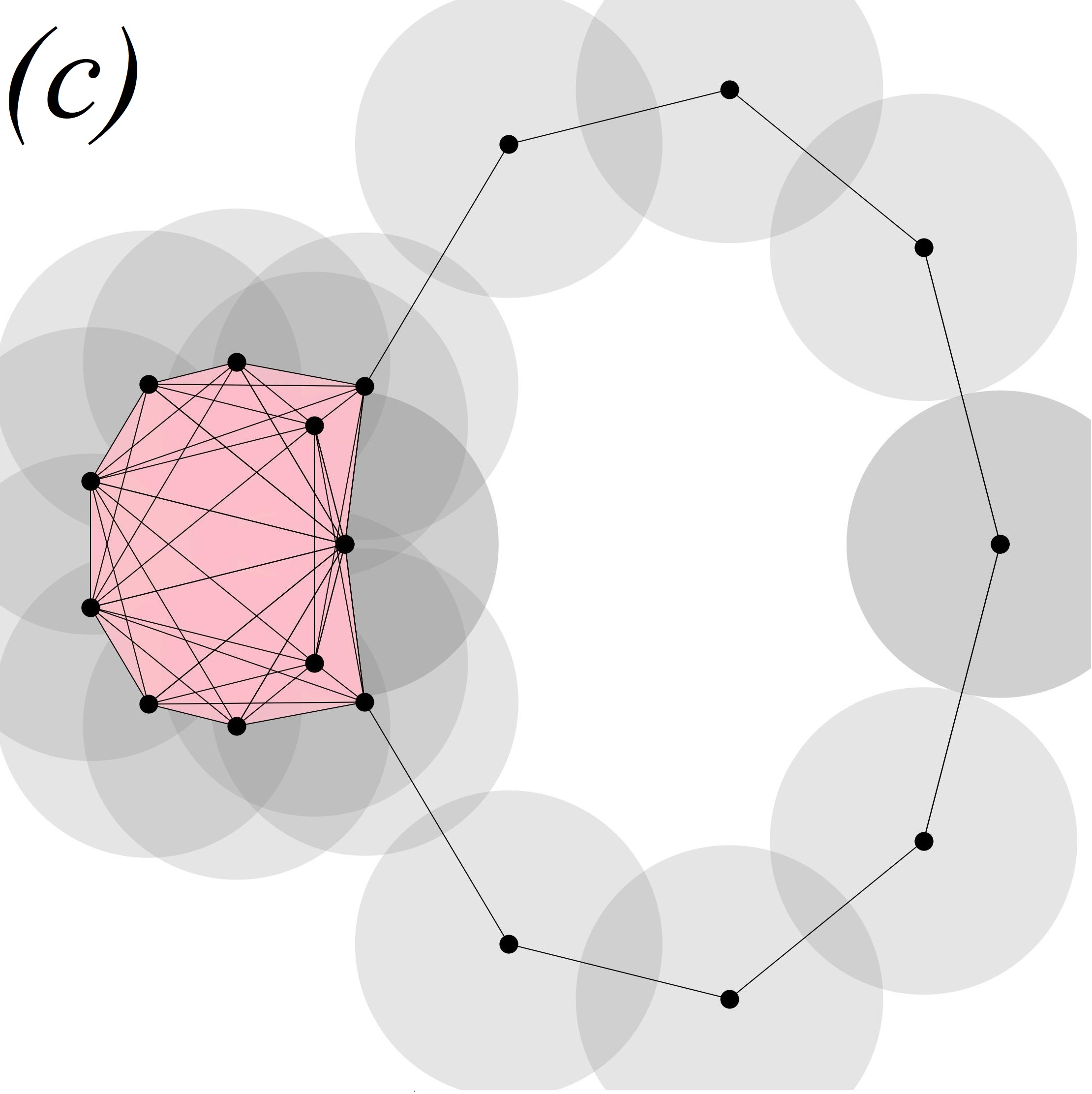}
    \includegraphics[height = 1in]{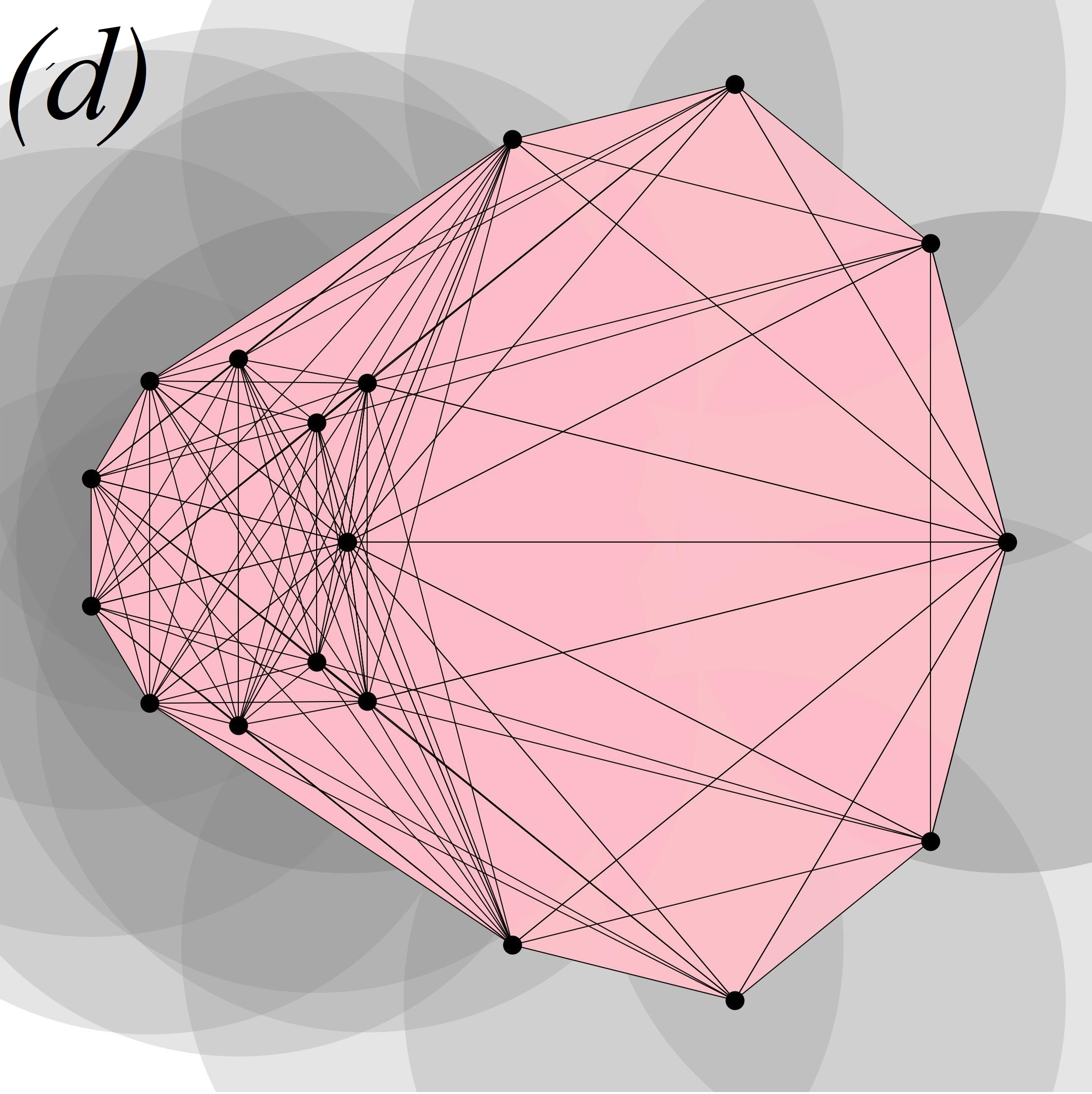}
    \caption{Sublevel persistence of a point cloud. Figures correspond to radii of (a) 0.14 (birth of smaller loop), (b) 0.31 (death of smaller loop), (c) 0.34 (birth of bigger loop), and (d) 0.73 (death of bigger loop).}
\label{fig:PointCloudPersistence}
\end{figure}

We will be particularly interested in sublevelset persistence in the case of a point cloud $\chi = \{x_1,\cdots,x_n\} \subset \R^d$ given as input. 
For this point cloud, and a fixed $r\geq 0$, the Rips complex is given by 
\begin{equation*}
    R(\chi,r) = \{ \sigma = \{ x_0,\cdots,x_d\} \mid \|x_i-x_j\| \leq r \text{ for all }i,j \}. 
\end{equation*}
That is to say, we treat the collection of points abstractly as a vertex set, and include a simplex for every subset of points that are pairwise within distance $r$ of each other. 
Note that if we can instead construct a function $f:K \to \R$ where $K$ is the complete simplicial complex on $n$ vertices (that is, the complex where every possible subset of the vertices is included), given by $f(\sigma) = \max_{x_i,x_j \in \sigma} \|x_i-x_j\|$. 
In this setting, the Rips complex is the same as the sublevelset $R(\chi,r) = f\inv(\infty,r]$. 
We thus construct a sublevelset filtration and the resulting persistence diagram accordingly. See the example in Fig.~\ref{fig:PointCloudPersistence}. As radius $r$ is varied, we get different Rips complexes, $R(\chi,r)$. Note the two points in the persistence diagram represent the two loops in the data. 

\subsection{Betti Number \& Vector}
Often times we are not actually interested in the full homology vector space, but instead are simply interested in the dimension. 
The dimension of homology is called the \textit{Betti number}, and is commonly denoted by $\beta_p(X) = \dim(H_p(X))$. 

For any filtered superlevel complex, $K^n \subseteq K^1 \subseteq \cdots K^0$, the Betti curve is a function defined as:
$L \mapsto \left(\beta_p (K^{\lfloor L \rfloor})\right)$
That is, the Betti curve is a function mapping to the $p$-th Betti number of the space at its filtration value of $L$. 
In the case of discretized data, we often have this data as a vector 
\[[\beta_p(L^0), \beta_p(L^1), \cdots, \beta_p(L^N)]\]
in which case, we refer to it as a Betti vector. In our case of cubical complexes, the filtration parameter $L$ corresponds to the function value (height in the PDF).

\subsection{Topological Consistency with Estimated Density}
\label{sec:topological_consistency}

In this work, we will be dealing with Kernel Density Estimates (KDEs) of system responses, and thus need to recover the homology of the analytical PDF from the KDE. For this, we follow the work of Bobrowski et al.~\cite{Bobrowski2017} to use the ideas of persistence in order to determine the homology of a single superlevelset of a KDE in the cubical complex case using a point cloud sampled from it.

Assume we are given a PDF $p:\R^d \to \R$, and let $D^L = p\inv([L,\infty))$ be a superlevelset of $p$. 
Our goal will be to recover the homology $H_p(D^L)$ given samples $\{X_1,\cdots,X_n\} \sim^{iid} p(x)$. 
Let $\hat p:\R^d \to \R$ be a kernel density estimate built using this sample (see \cite{Scott1992-ld} for multivariate density estimation). 
We then estimate the set $D^L$ using the points from the sample by using a union of balls $B_r(x) = \{y \in \R^d \mid \|x-y\| <r\}$ centered at points in the original which have estimated value above $L$. 
Specifically, let $\chi ^{L} := \{X_{i}: \hat{p}(X_{i}) \geq L\}$, and then this approximation to $D^L$ is given by 
\[
\widehat{D}^L_r := \bigcup_{X \in \chi^L}  B_{r}(X).
\]
Following~\cite{Bobrowski2017}, we choose the radius $r$ to be the same as the bandwidth used for computing the kernel estimator. 

Given $\e>0$, the inclusion map
$ i: \hat{D}^{L+\e}_r \xhookrightarrow{} \hat{D}^{L-\e}_r $
induces a map on the homology
\[i_*: H_{p}(\hat{D}_{L+\e}(r)) \xrightarrow{} H_{p}(\hat{D}_{L-\e}(r)).\]
This map can be used to define an estimator for the homology of superlevel-set at $L$ by setting
\[\widehat{H_{*}}(L, r):=\im(i_{*}).\]
In Bobrowski et al.~\cite{Bobrowski2017}, it is shown that for a large enough sample and good choices of $\e$ and $r$, $\widehat {H_p}(L,r) \iso H_p(D^L)$ with high probability. The details of the algorithm can be found in their paper~\cite{Bobrowski2017}. Although there is no definitive method for estimating the minimum sample size and \textit{good choices} for $\e$ and $r$, the following guidelines outlined in~\cite{Bobrowski2017} should be kept in mind: 

\begin{enumerate}
    \item While the algorithm does not specify a lower bound for the parameter $\e$, it must be chosen so that there are no critical values in the range $[L-2\e,L+2\e]$. This, of course, is not quite possible in practice and will be inadvertently violated by our method. Therefore, we expect to see errors occur near critical points where the homology changes. Regardless, these errors do not damage the method's capability to detect changes in system PDF and the onset of bifurcation, as you will see later.
    
    \item For the parameter $r$, the user must be cautious of the lower bound on it, given by $nr^{d} >> log(n)$ (in stricter terms). However, not all values above this lower bound will estimate the correct homology, hence caution needs to be taken while applying the algorithm.
\end{enumerate}

For our examples, we have kept $\e$ in the order of $10^{-5}$, and $r$ between $0.1$ and $0.8$ (varies case to case), against an $n = 500$ ($500$ sparse values extracted from a larger set using greedy permutation \cite{greedyPerm}). 

\section{Methods}

In Sec.~\ref{sec:maths}, we discussed how a single persistence diagram from one cubical complex can be distilled to its betti vector. For a family of cubical complexes changing over an external variable such as the bifurcation parameter, a sequence of persistence diagrams can be realized. The evolution of such persistence diagrams can be visualized in the bifurcation interval by plotting the changing Betti vectors against the bifurcation parameter. We shall refer to these plots---analogous to the CROCKER plots devised by Topaz et al.~\cite{CROCKER} to visualize the evolution of sublevel filtrations of Rips complexes---as the "homological bifurcation plots'', and see how they help us detect a P-bifurcation. 

To explain and illustrate the methodology, we rely on the example of the stochastic Duffing oscillator. Analytical expressions for the stationary joint PDF for the oscillator was derived by splitting the FPK equation by Mamis and Athanassoulis \cite{Mamis2016}.

The stochastic Duffing Oscillator forced by an additive white gaussian noise is represented by the SDE
\[\ddot{X} + \dot{X} + hX + X^{3} = q_{1}dW_{1}.\]
The stationary joint PDF for the oscillator is given by:
\[
p_{x_{1}x_{2}}(\textbf{x}) 
= C\exp\left[-\frac{1}{2q^{2}_{1}D_{11}}
    \left(x^{2}_{2} + hx^{2}_{1} + \frac{1}{2}x^{4}_{1}\right)
\right]
\]
where $C$ is a normalizing constant such that $\int_{-\infty}^\infty pdx_{1}dx_{2} = 1$; and $D_{11}$ is from the intensity matrix $D_{ij}$ of the Gaussian white noise and is proportional to its power spectral density. Since the effect of changing noise is not being investigated here,  we have kept the noise amplitude and intensity constant: $q_{1} = D_{11} = 1$. For ease of computation, we have used $C = p/$max$(p)$---which does not affect the topology of the PDF being a scalar multiplier.

\begin{figure}[!htbp]
    \centering
    \includegraphics[width = 0.19\textwidth]{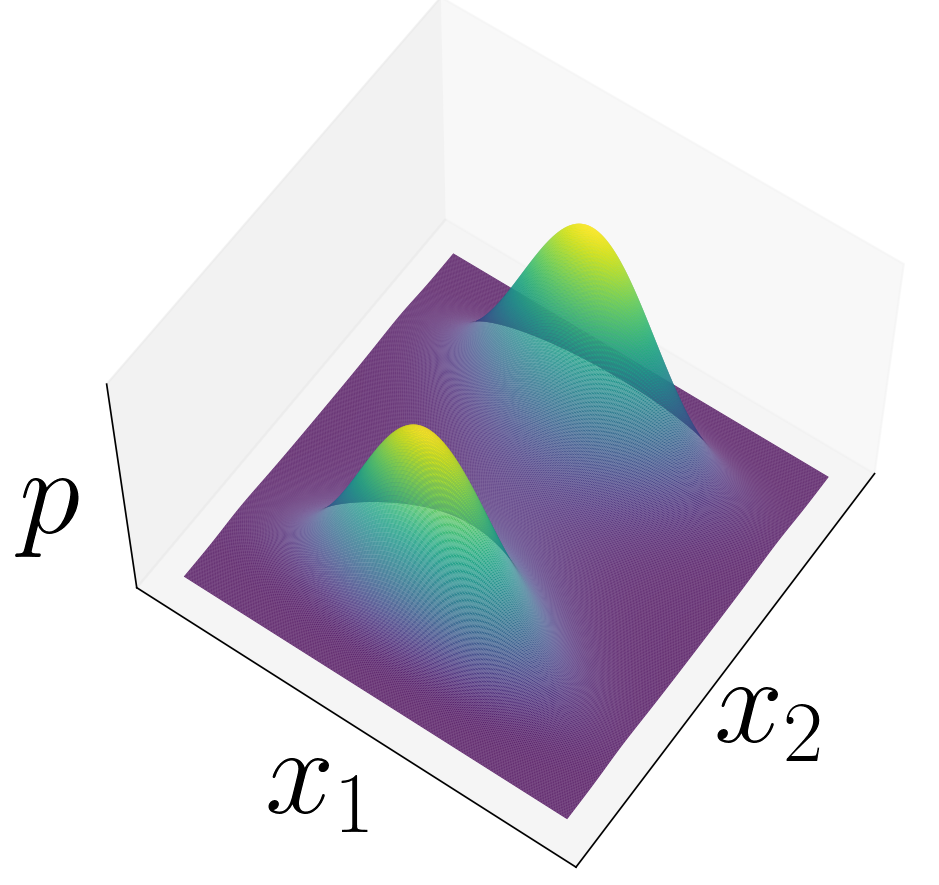}
    \includegraphics[width = 0.19\textwidth]{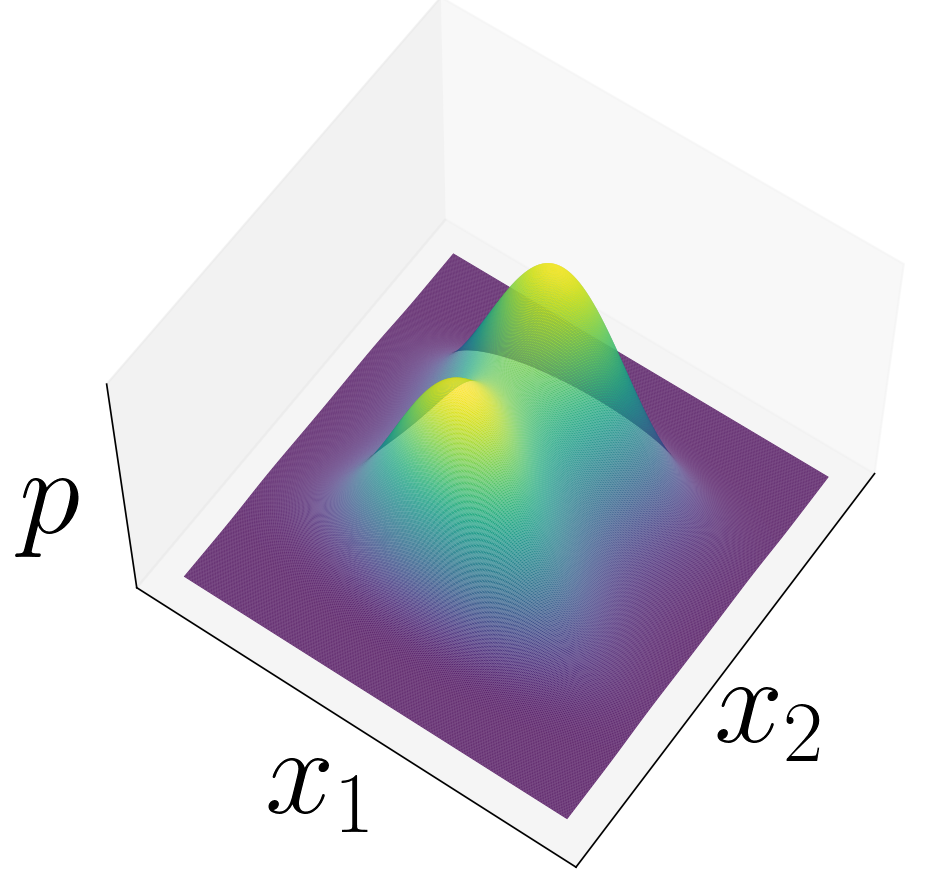}
    \includegraphics[width = 0.19\textwidth]{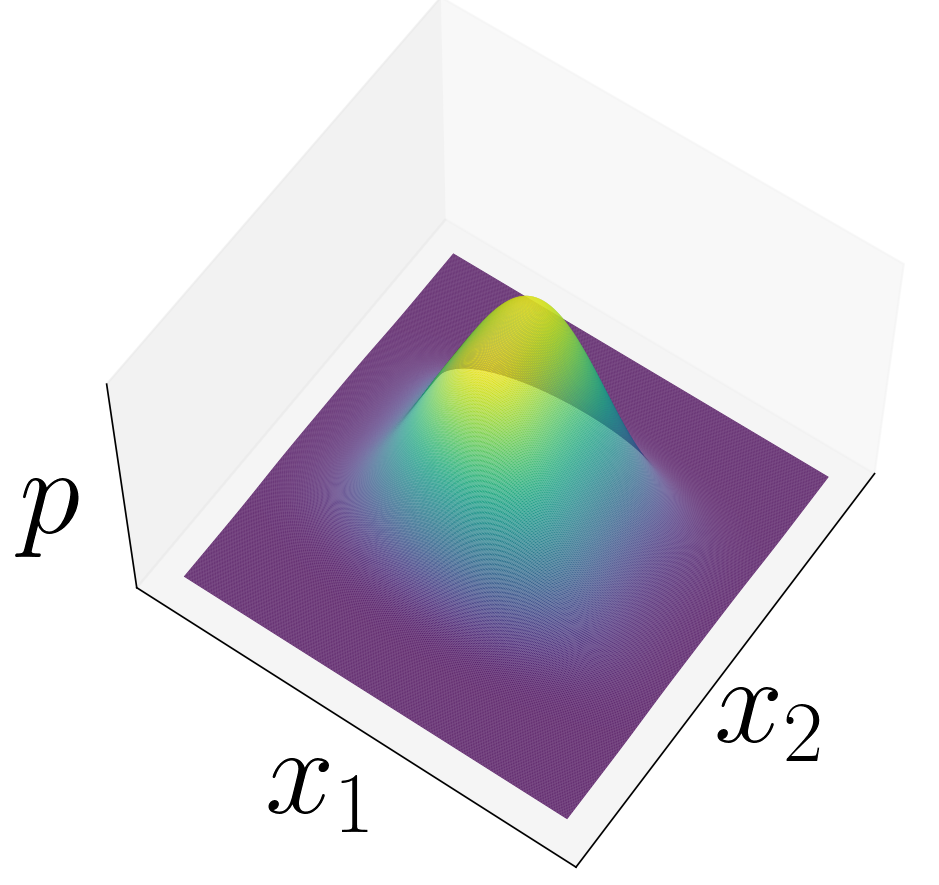}
    \includegraphics[width = 0.19\textwidth]{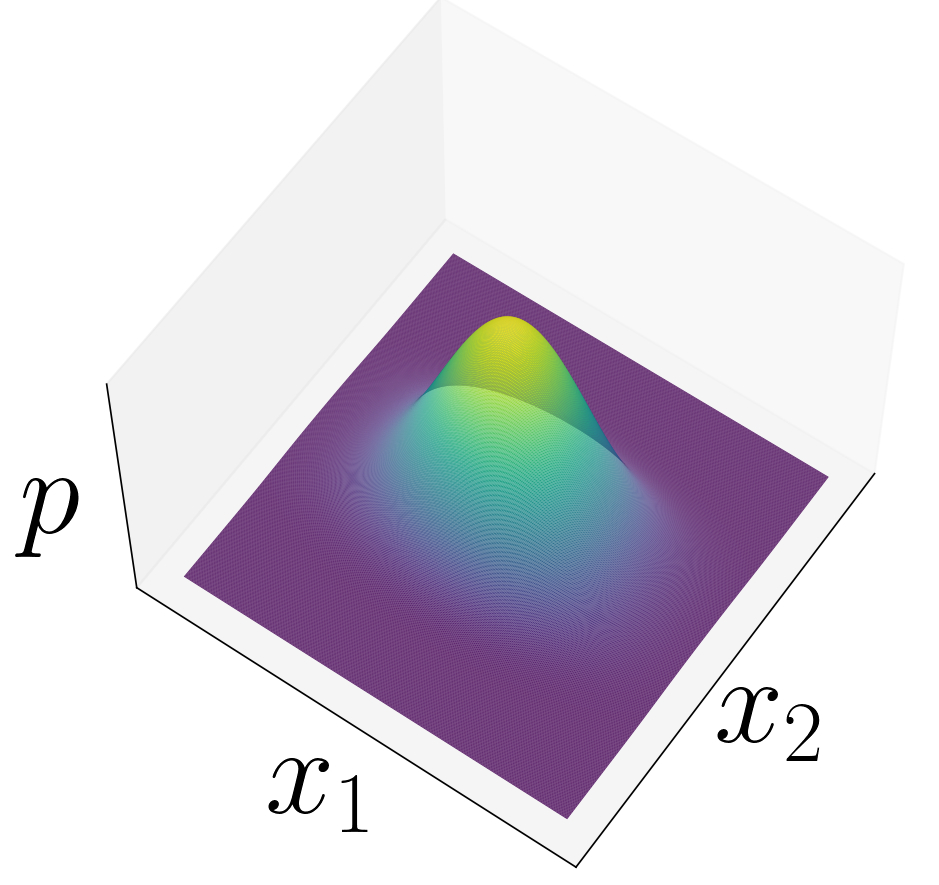}
    \includegraphics[width = 0.19\textwidth]{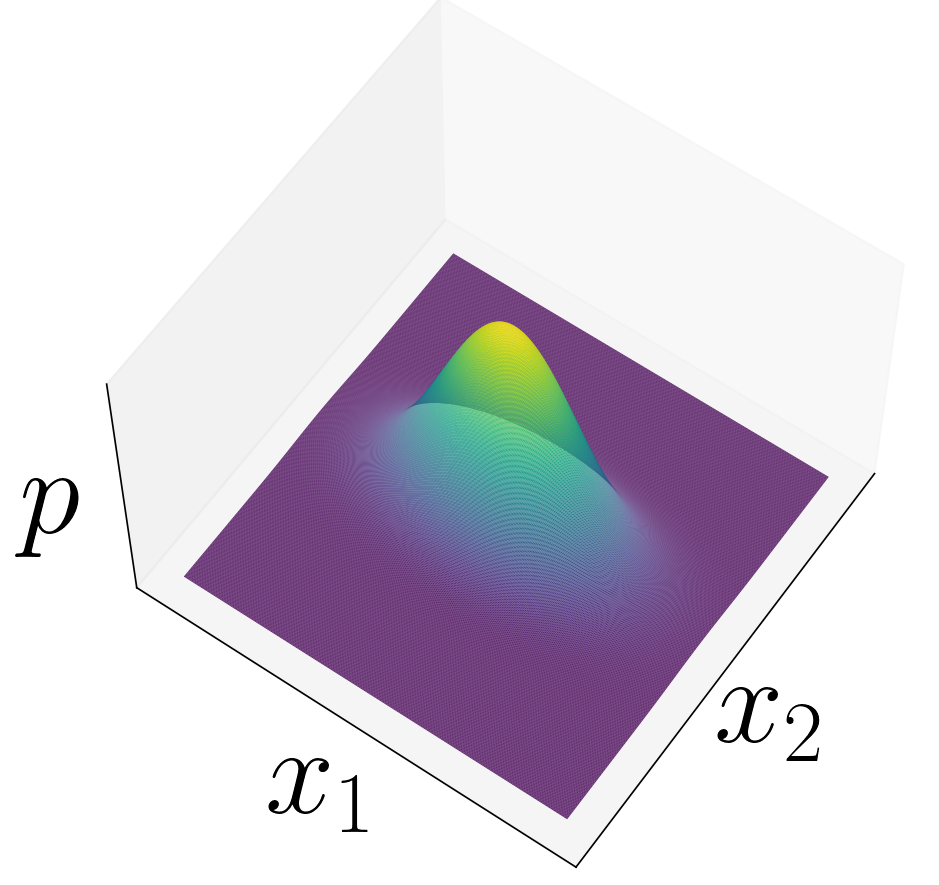}
    \caption{Evolution of the PDF for stochastic Duffing oscillator as $h$ changes its values. From left to right, $h$ = $-3, -1, 0, 1, 3$.}
        \label{fig:Duffing_analytical}
\end{figure}

For this oscillator, as the bifurcation parameter $h$ goes from negative to positive, we notice the shape of the PDF changing from bi-stability to mono-stability. Figure~\ref{fig:Duffing_analytical} shows the evolution of the analytically generated PDFs corresponding to a variety of values for the bifurcation parameter, $h$. Since the purpose of detecting a bifurcation would still be satisfied between $-1$ and $1$, the next sections only deal with $h \in [-1, 1]$.
\subsection{Detection using Analytical PDF}
\begin{figure}[!htbp]
    \centering
    \includegraphics[width = 0.2\textwidth]{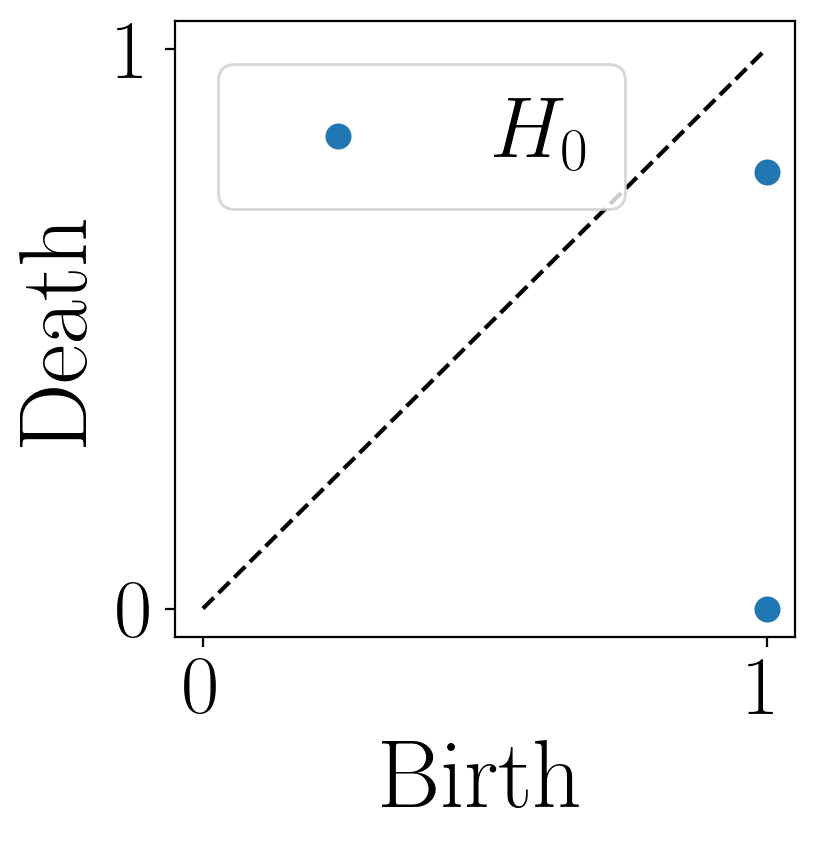}
    \hfil
    \includegraphics[width = 0.2\textwidth]{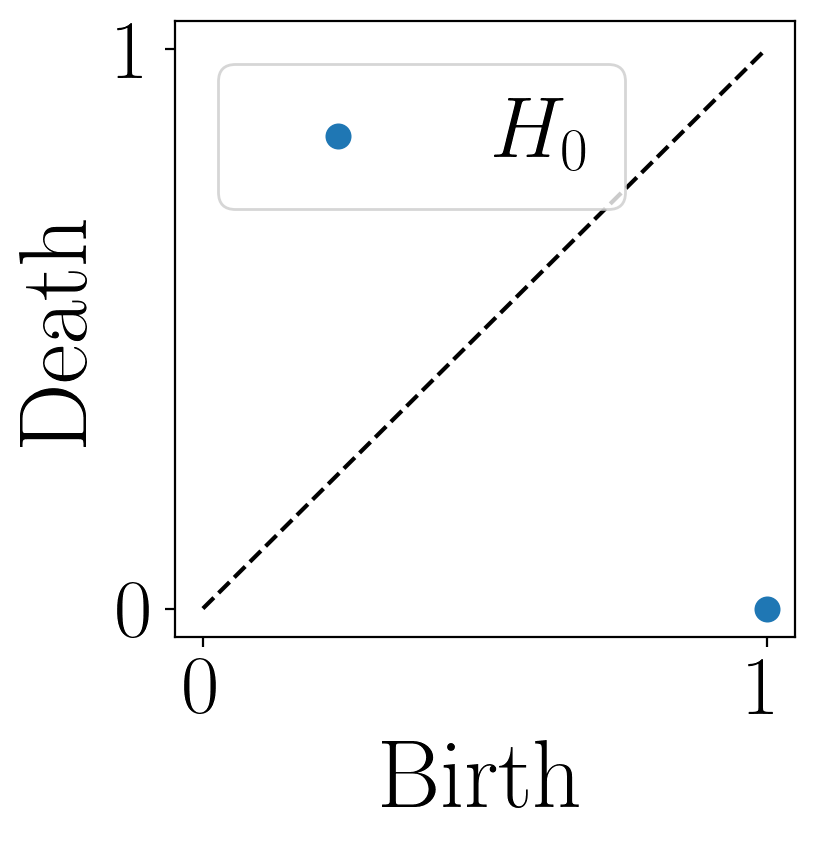}
    \hfil
    \includegraphics[width = 0.2\textwidth]{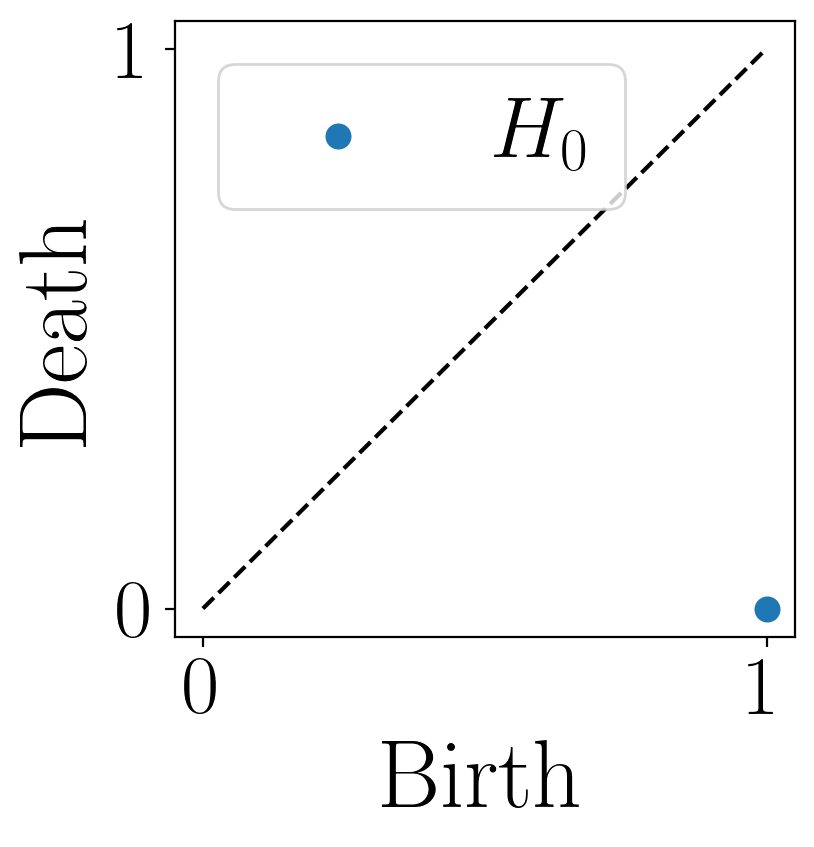}
    \caption{Effect of bifurcation on the rank of $H_{0}$ for Duffing oscillator. From left to right, $h = -1, 0, 1$.}
    \label{fig:Duffing_analytical_bif}
\end{figure}

Given the analytical PDF for Duffing oscillator, we compute the super-level persistence diagrams for the cubical complexes corresponding to $h \in \{-1, 0, 1\}$. Notice in Fig.~\ref{fig:Duffing_analytical_bif}, the rank of $H_{0}$ changes from $2$ to $1$ at $h = 0$ (when the bifurcation occurs) and stays $1$ afterwards. Each $H_{0}$ component corresponds to a peak ("connected component") in the PDF. 

\begin{figure}[!htbp]
\centering
\includegraphics[width=0.7\linewidth]{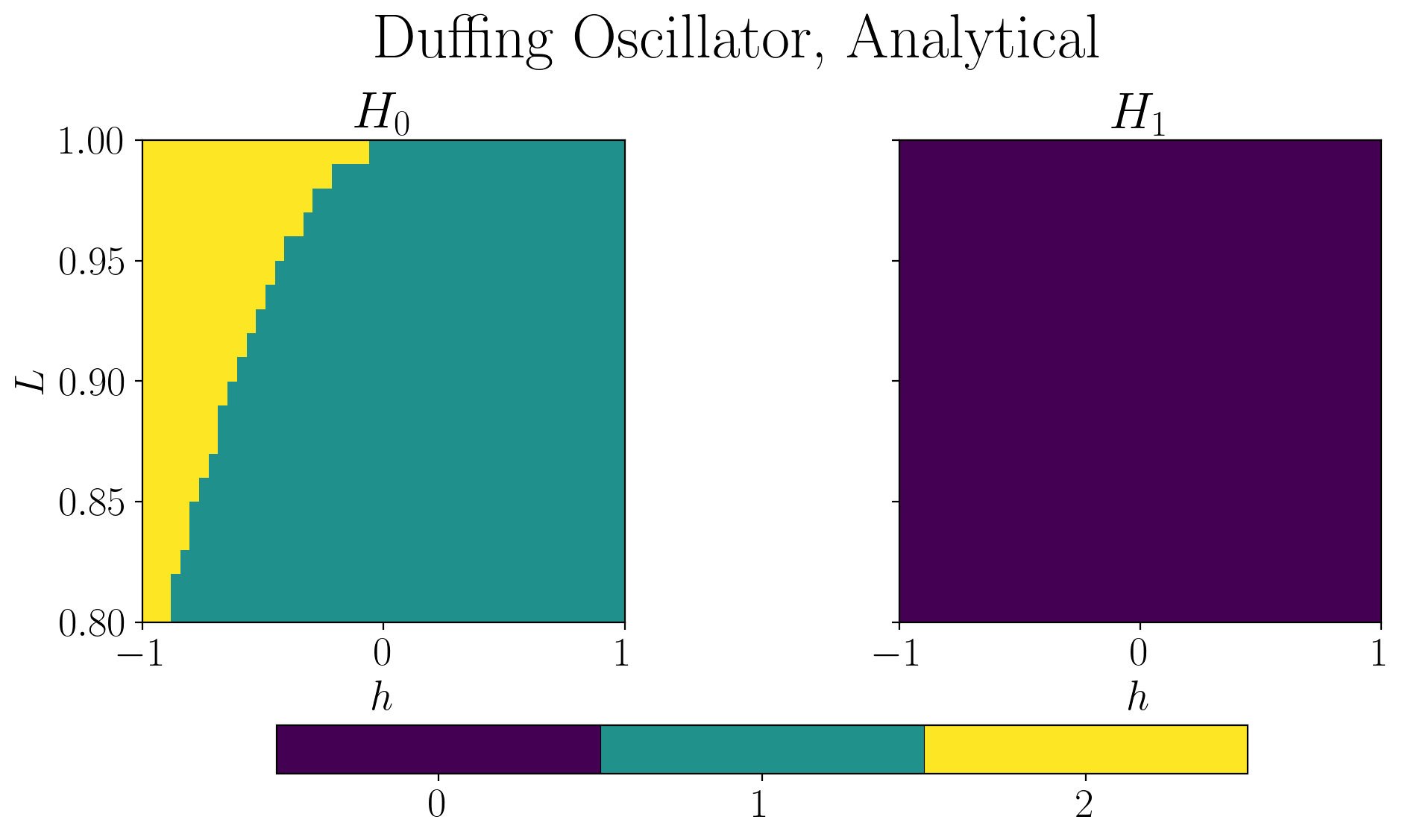}
\caption{Analytical homological bifurcation plots plots for Duffing Oscillator.}
\label{fig:duffing_analytical_crocker}
\end{figure}

Thus, bifurcations can be linked to the change in homology rank for the associated complex. Figure~\ref{fig:duffing_analytical_crocker} shows the homological bifurcation plots generated from the analytical PDFs for Duffing Oscillator over the interval $h \in [-1, 1]$. The $x$-axis of the plot corresponds to the bifurcation parameter $h$, the $y$-axis corresponds to the level $L$ and the colour of the points dictates the Betti number at that level. Notice the change in ranks initiating at the maximum level when the parameter turns non-negative. Unlike this example, in a case where the PDF bifurcates into a Limit Cycle, the change would occur in the rank of $H_{1}$ group (as would be shown later).

\subsection{Detection using Estimated Density}
As opposed to the previous section, in real-world scenarios, the analytical PDF of the system response is rarely available. However, realizations of the state vector of the dynamical system $X = [x_1,\cdots,x_N]$ where $x_i \in \R^p$ ($p$ being the components in state space) may be sampled. With such an input, the system response can be turned into a KDE \cite{Scott1992-ld}. %
However, neither acquired system data nor the KDEs are above being noisy, and given the discrete nature of homology, a single inaccurately placed noise artefact may change the homology rank, as demonstrated here \cite{Bobrowski2017}. This is where the discussion in Sec.~\ref{sec:topological_consistency} comes in for computing an estimate of the betti number at a specific level of the KDE.

\begin{figure}[!htbp]
    \centering
    \includegraphics[width = 0.7\textwidth]{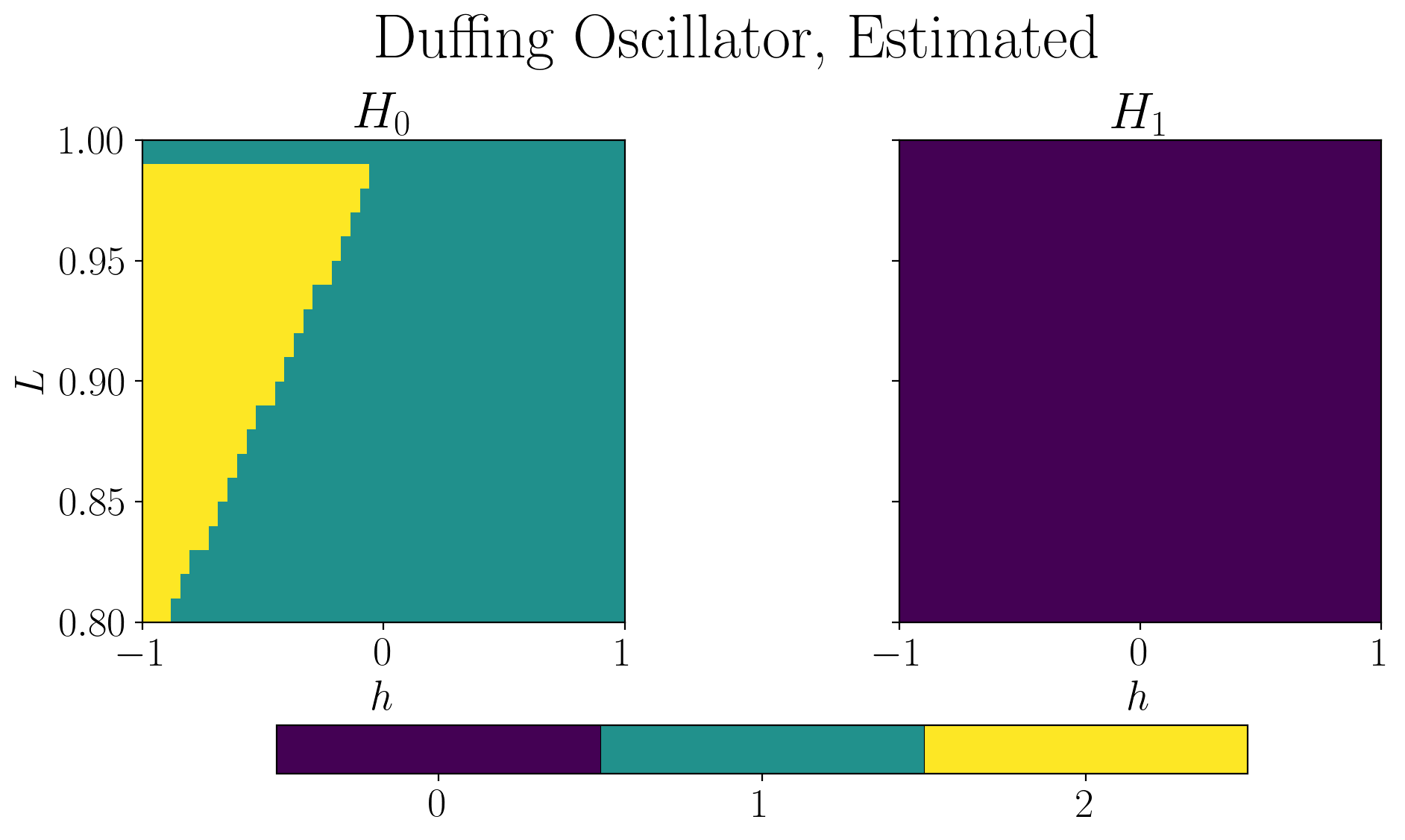}
    \caption{Estimated homological bifurcation plots for Duffing oscillator.}
\label{fig:Duffing_short_CROCKER}
\end{figure}

Once we have access to a KDE against one value of the bifurcation parameter, the homology ranks can be computed for a range of levels, $L$, and plotted to visualize the change in ranks. 
Figure~\ref{fig:Duffing_short_CROCKER} shows the estimated homological bifurcation plots for the betti numbers for $H_0$ and $H_1$ for Duffing over a range of values for the bifurcation parameter.
The estimated plot matches the analytical and indicates a bifurcation around $h = 0$ in $H_0$. As expected, there is no change in the plot for $H_{1}$.

\subsubsection{Error between Analytical and Estimated Homological Bifurcation Plots}

\begin{figure}[h]
\centering
\includegraphics[width=0.7\linewidth]{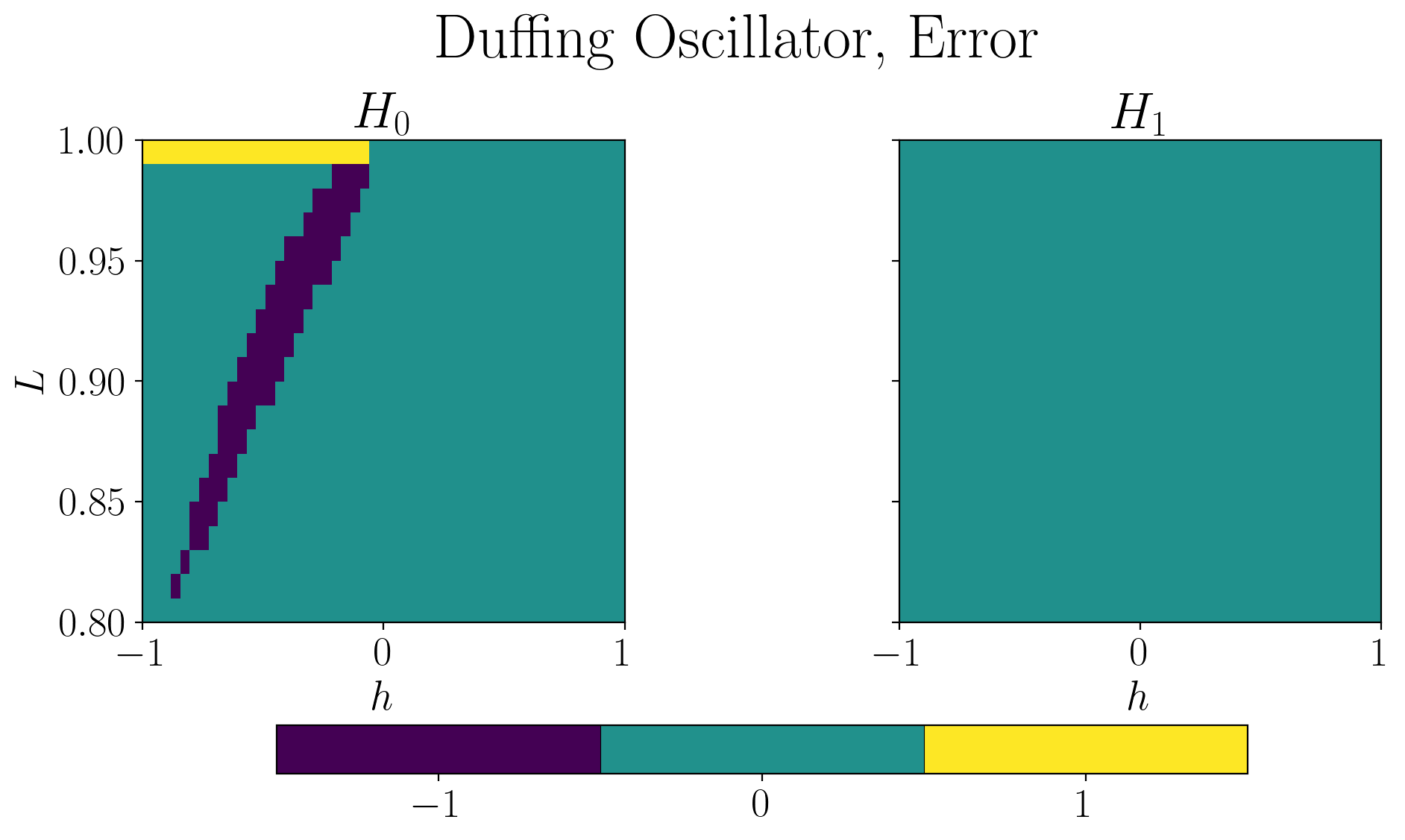}
\caption{Error between analytical and estimated homological bifurcation plots for Duffing.}
\label{fig:duffing_error_crocker}
\end{figure}

Figure~\ref{fig:duffing_error_crocker} shows the error between the analytical and estimated homological bifurcation plots for Duffing. Since homology is discrete, the error is computed as $\beta_\text{true} - \beta_\text{estimate}$. Note that the analytical and estimated values do not match around a `topological boundary' in the plots. All errors (at $L \approx 1$ and the triangular region), are introduced due to those levels corresponding to critical points (topological boundaries) in the PDF. %

\section{Results and Discussion}
\label{sec:results}

In this section, we demonstrate the successful application of the methodology on two other stochastic oscillators and discuss the results.

\subsection{Raleigh-Vander Pol (RVP) Oscillator}
\begin{figure}[!htbp]
    \centering
    \includegraphics[width = 0.23\textwidth]{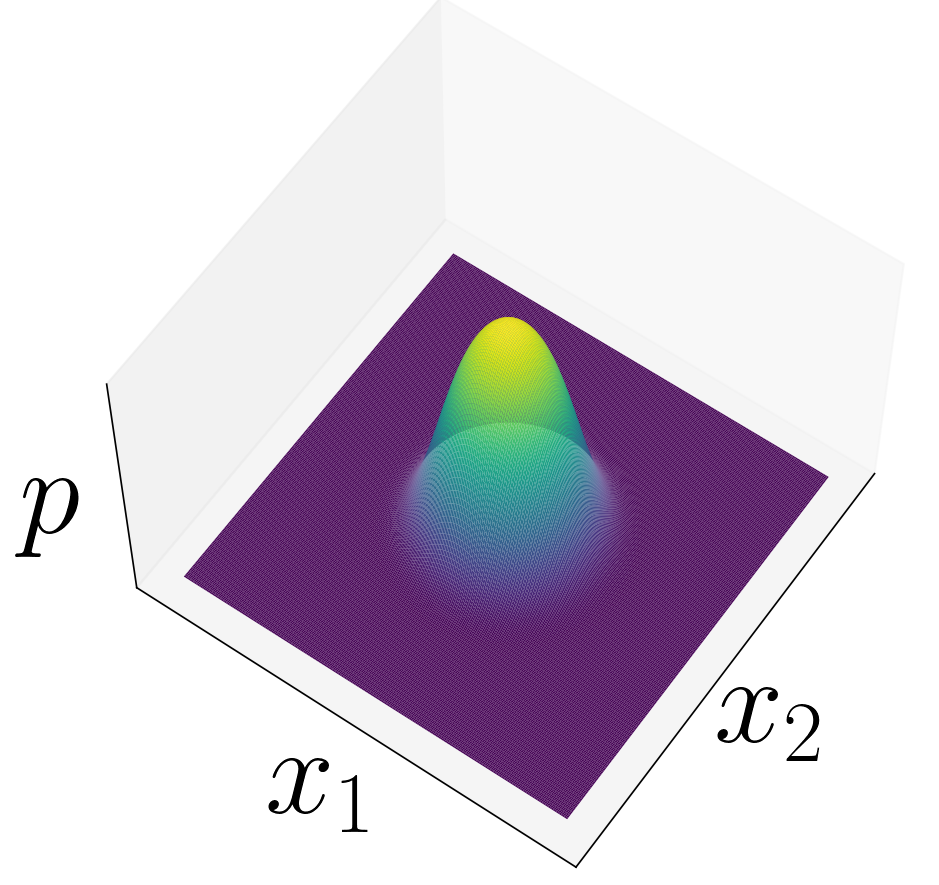}
    \hfil
    \includegraphics[width = 0.2\textwidth]{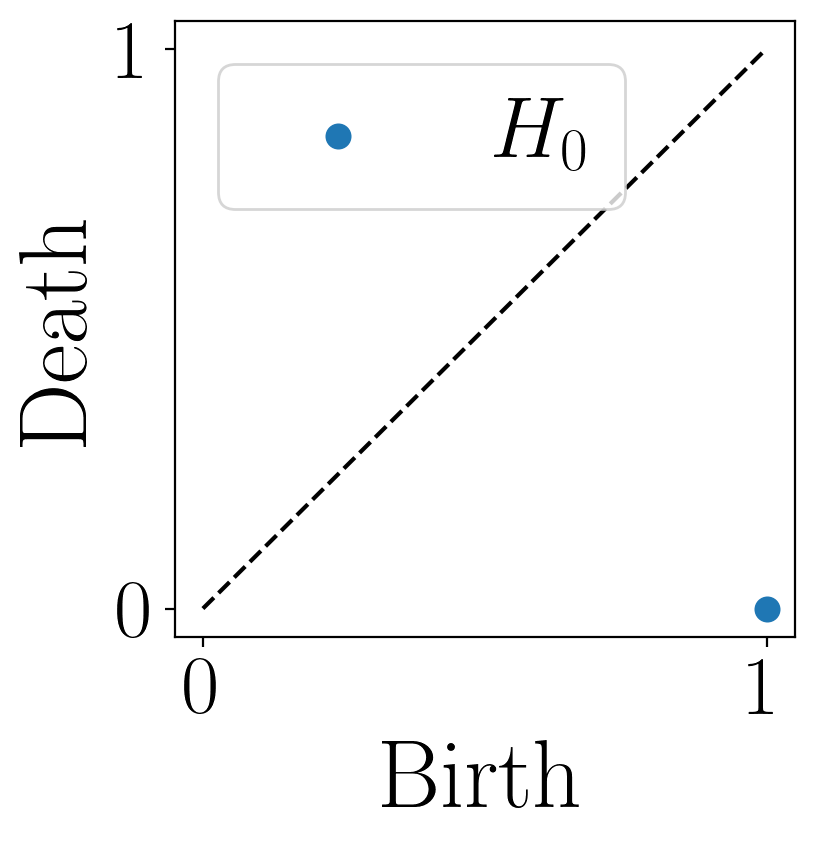}
    \hfil
    \includegraphics[width = 0.23\textwidth]{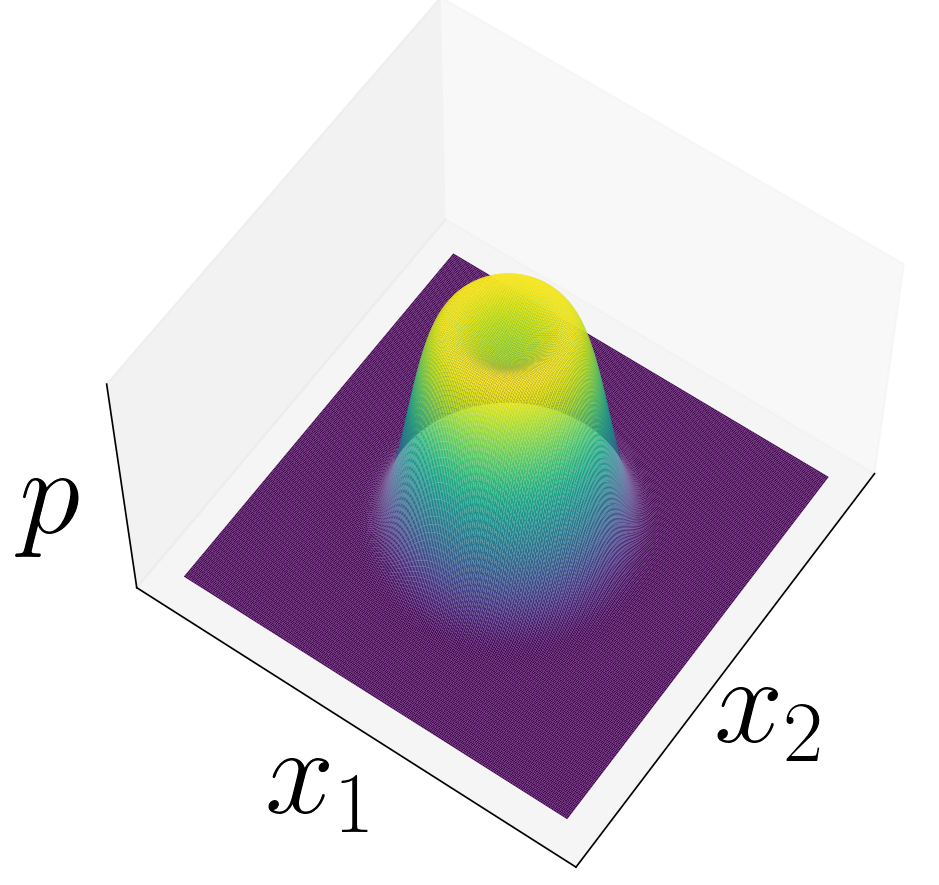}
    \hfil
    \includegraphics[width = 0.2\textwidth]{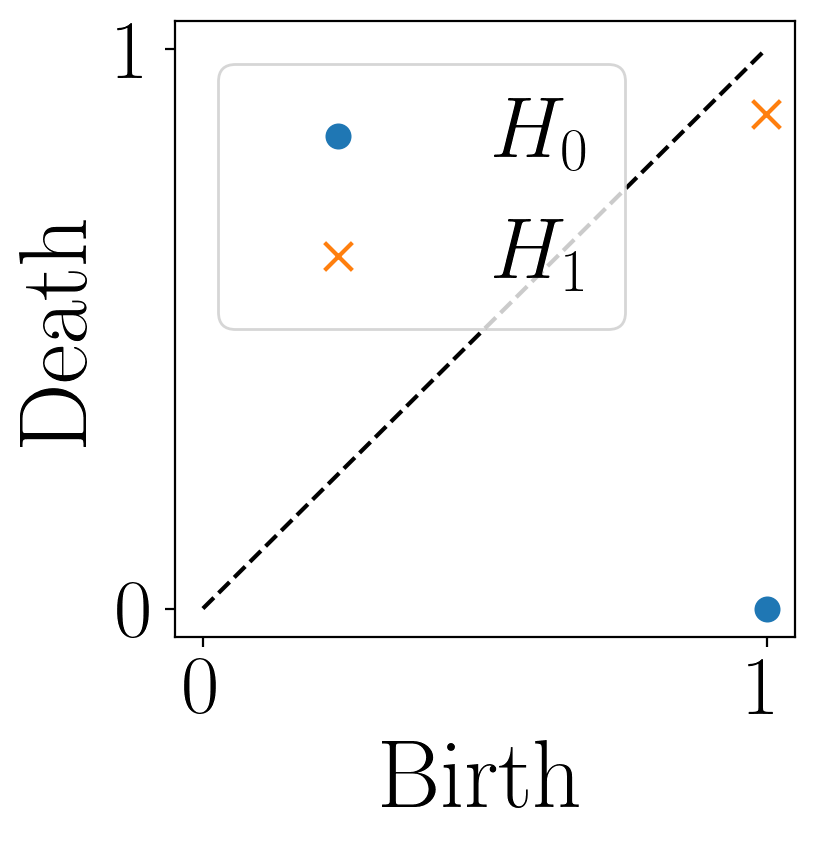}
    \caption{
        PDFs for Monostable ($h=1$, left) and Limit Cycle ($h=-1$, right) RVP oscillator respectively, with their superlevel cubical persistence diagrams.}
        \label{fig:Vander_analytical}
\end{figure}

A stochastic Rayleigh-Vander Pol (RVP) Oscillator forced by an additive white gaussian noise is represented by the following SDE:
\[\ddot{X} + (h + X^{2} + \dot{X}^{2})\dot{X} + X = q_{1}dW_{1}\]
where C and $D_{11}$ are the same as defined earlier. 
The stationary joint PDF \cite{Mamis2016} is as follows:
\[
p_{x_{1}x_{2}}(\mathbf{x}) 
= C \exp \left[
-\frac{1}{2q^{2}_{1}D_{11}}
\left\{\frac{1}{2}(x^{2}_{2}+x^{2}_{1})^{2} + h(x^{2}_{2}+x^{2}_{1})\right\}
\right]
\]
This oscillator experiences a bifurcation at $h=0$ when the system shifts from a monostable state to limit cycle oscillations (LCO) when $h<0$. Figure~\ref{fig:Vander_analytical} shows the analytically generated PDFs corresponding to $h = 1$ and $h = -1$ respectively, along with their super-level cubical persistence diagrams. The corresponding persistence diagrams show one $H_0$ component against one peak and one $H_1$ component corresponding to a crater or ``hole" in the PDF.

\begin{figure}[!htbp]
\centering
\includegraphics[width=0.6\linewidth]{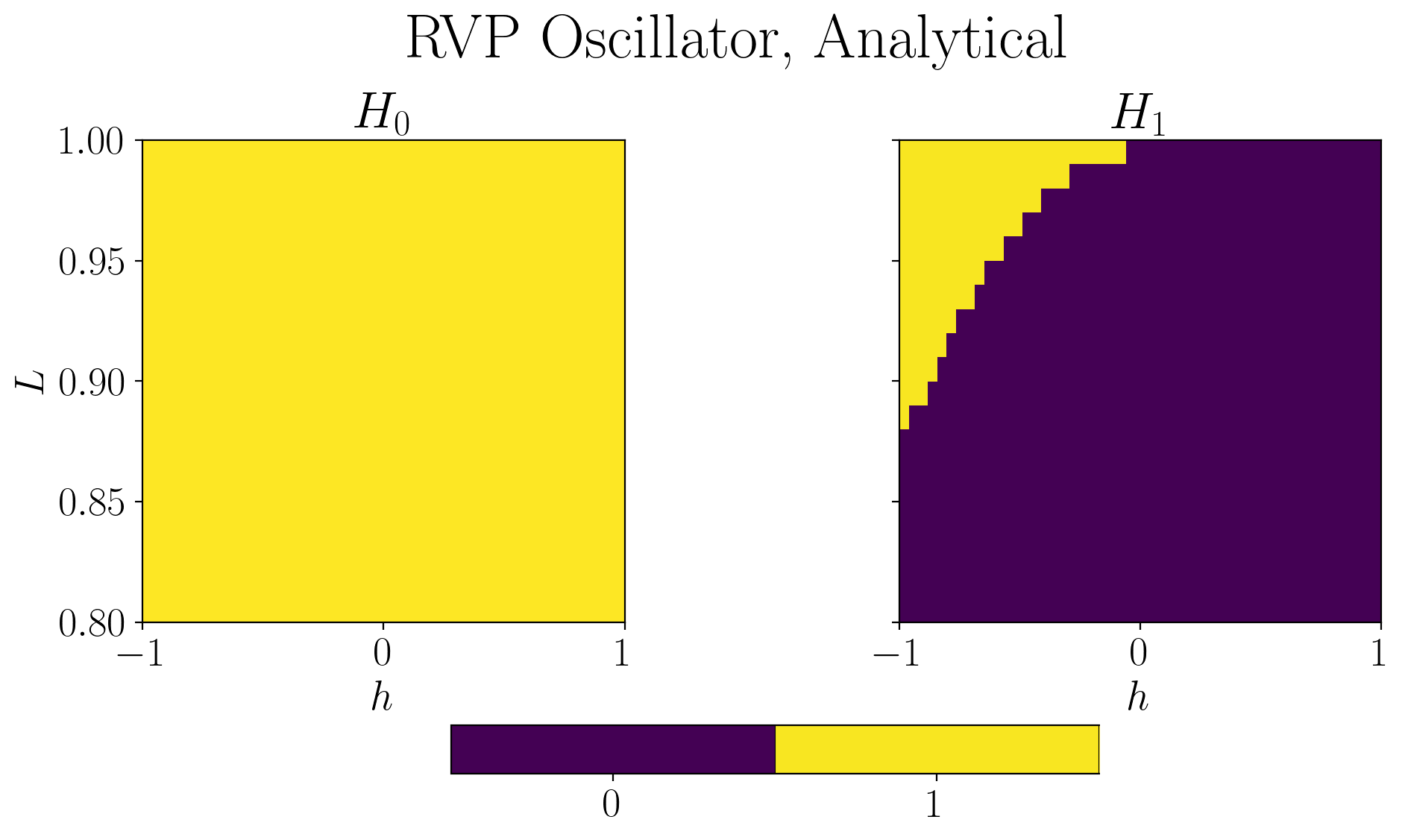}
\includegraphics[width = 0.6\textwidth]{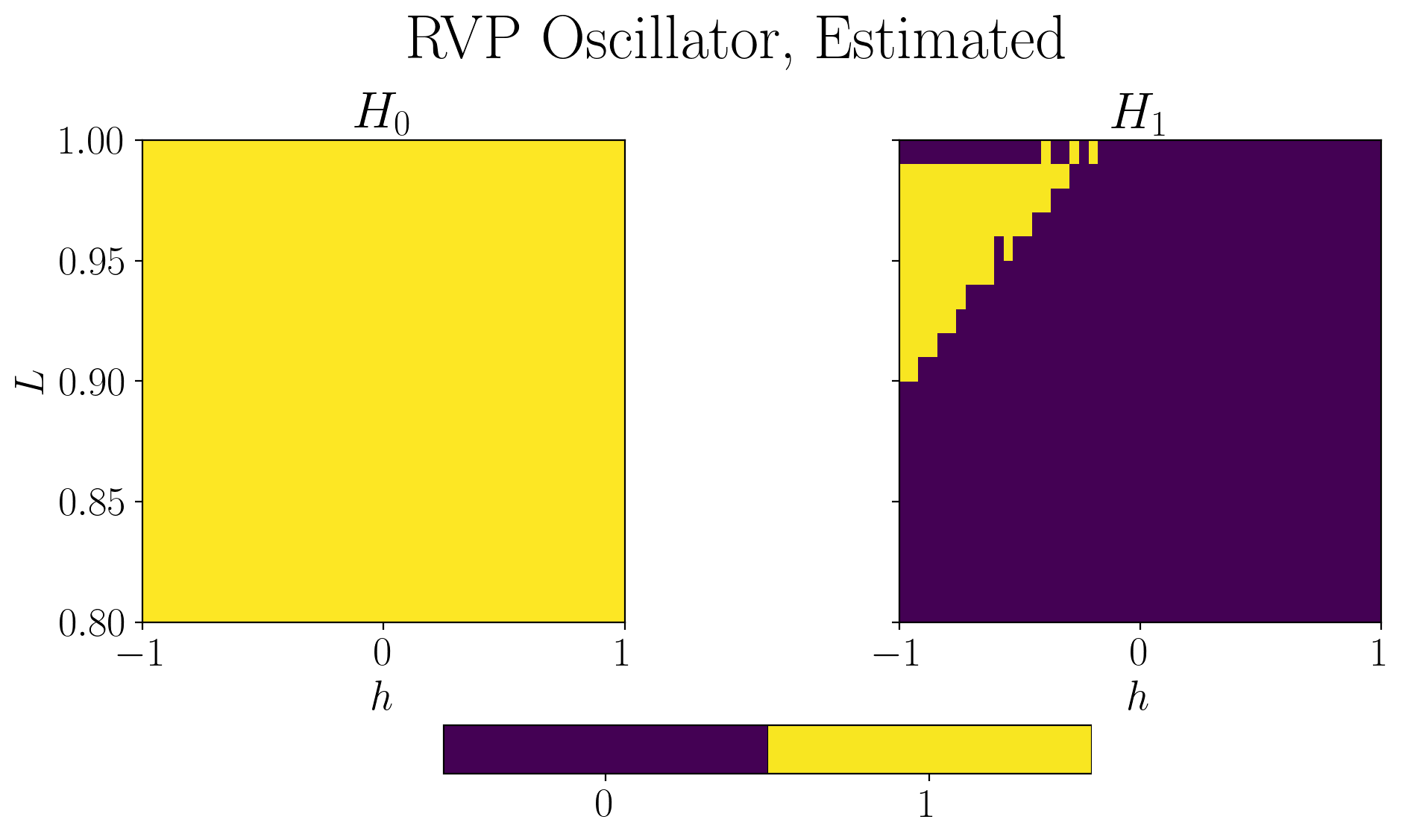}
\includegraphics[width=0.6\linewidth]{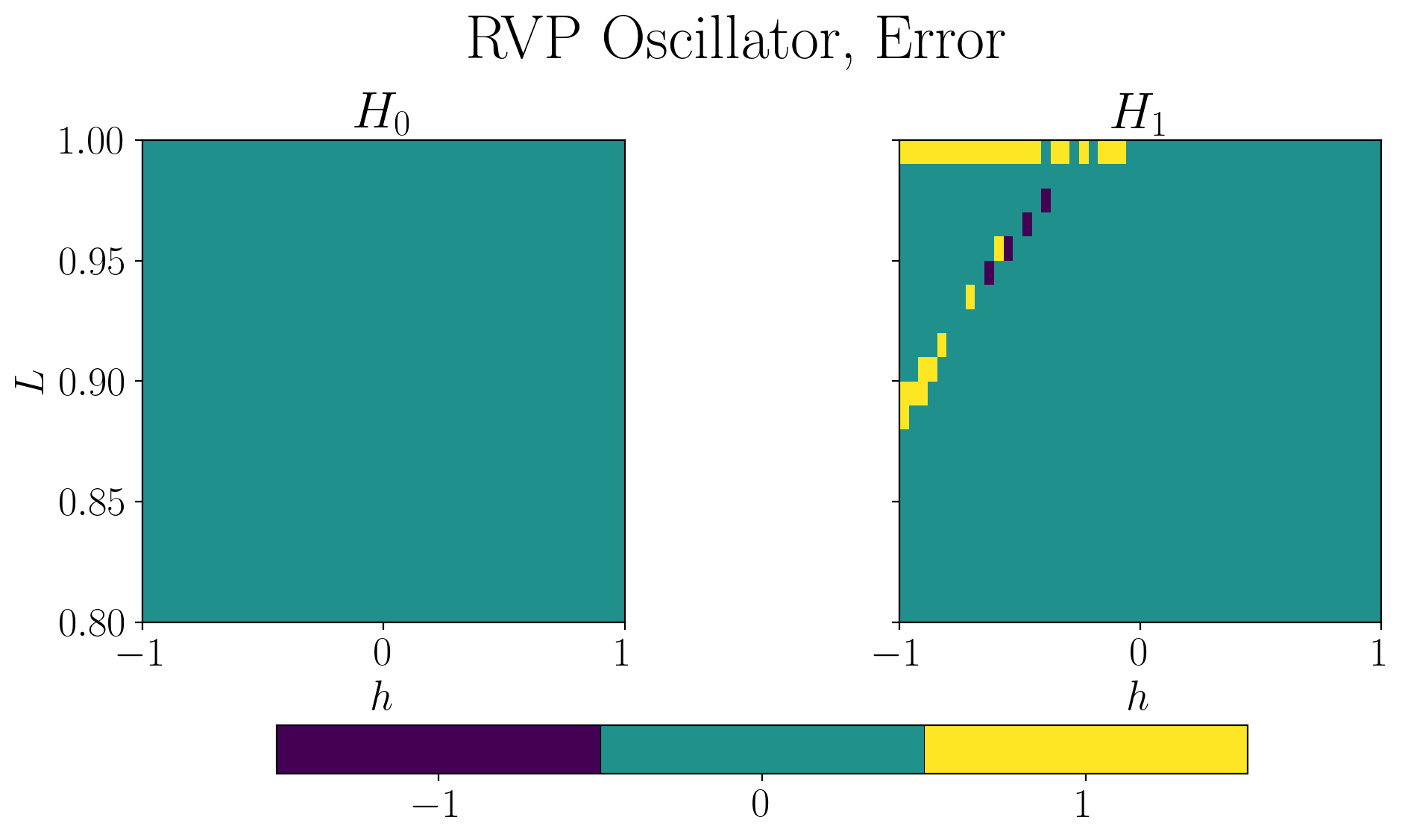}
\caption{Analytical and estimated homological bifurcation plots for RVP Oscillator, along with the errors.}
\label{fig:vander_analytical_crocker}
\end{figure}

Figure~\ref{fig:vander_analytical_crocker} shows the analytically generated and estimated homological bifurcation plots for RVP Oscillator over the interval $h \in [-1, 1]$, along with the errors between the two. Notice the change in homology when the parameter turns non-negative. In each vertical Betti vector, the $L$ at which the homology changes would correspond to the depth of the crater in the PDF. Estimated plot for $H_1$ shows similar features as the analytical plot and indicates the occurrence of a bifurcation. Finally, note that the errors between the analytical and estimated homological bifurcation plots occur at levels where the `topological' boundary is expected in the PDF. Similar to Duffing, at $L \approx 1$, some error is introduced due to the set $\hat{D}_{L+\e}$ being empty or nearly empty. %

\subsection{Quintic Oscillator}
Finally, we turn our attention to a slightly more complex stochastic oscillator called Quintic, which is forced by both an additive and a multiplicative white Gaussian noise. The oscillator is modelled by the following:
\[\ddot{X} + 2[D_{11}(2U(X) + h) - \frac{1}{2}D_{22}]\dot{X} + 2[D_{22}(2U(X) + h)\]
\[+ D_{11}]\dot{X}^{3} + 2D_{22}\dot{X}^{3} + h_{0}(X) = dW_{1} + \dot{X}dW_{2}\]
where $h_{0}(X) = x^{3}_{1} + ax^{2}_{1} - x_{1}$.
The stationary joint PDF for the system \cite{Mamis2016} is given by
\[p_{x_{1}x_{2}}(\textbf{x}) = C\text{exp}[-\frac{1}{2}x^{4}_{2} - (h + 2U(x_{1}))(x^{2}_{2} + U(x_{1})) - hU(x_{1})]\]
where C is the same as earlier, and $D_{11}$ and $D_{22}$ are from the intensity matrix $D_{ij}$ of the Gaussian white noise, while %
$U(x_{1}) = \frac{1}{4}x^{4}_{1} + \frac{a}{3}x^{3}_{1} - \frac{1}{2}x^{2}_{1}$ is the potential energy. 

This system has two bifurcation parameters ($h$ and $a$), and the form of the PDF can be complicated, exhibiting both (symmetric and asymmetric) bi-stable and limit-cycle behaviours. The symmetry of potential energy results in symmetric shapes for PDFs, while negative values in damping parameters result in pronounced limit cycle behaviour. 
\begin{figure}[!htbp]
    \centering
    \includegraphics[width = 0.45\textwidth]{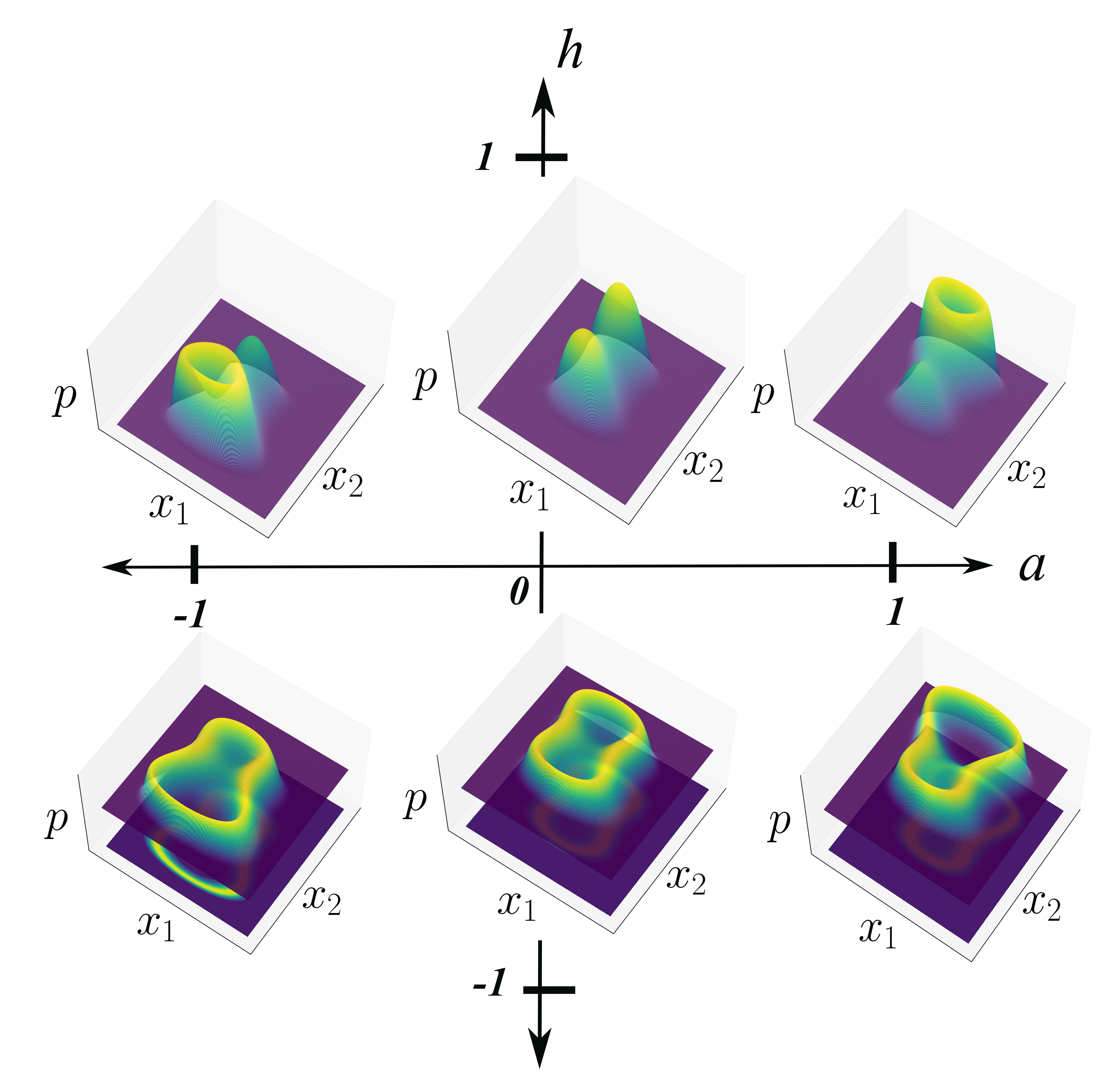}
    \includegraphics[width = 0.45\textwidth]{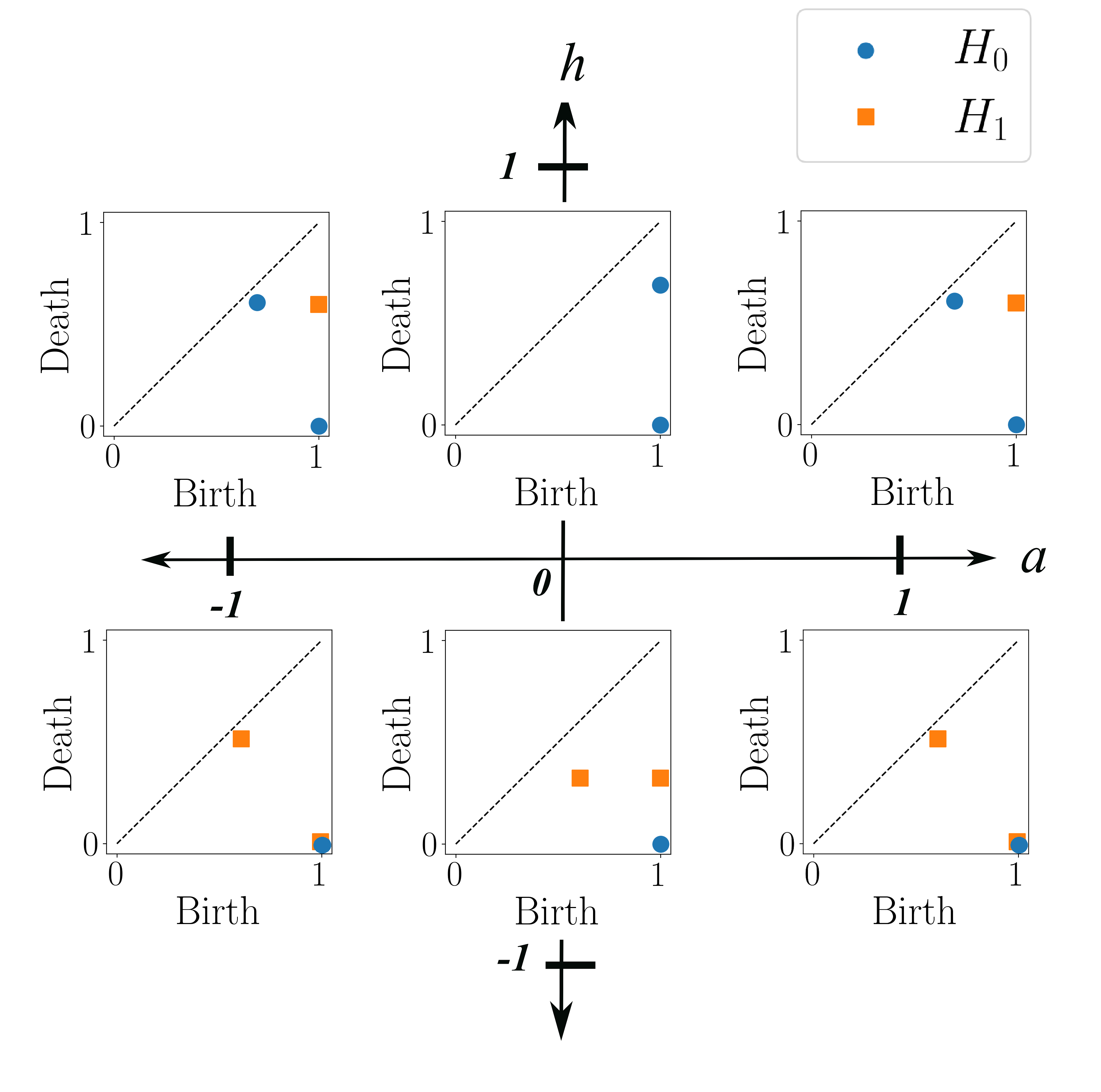}
    \caption{
        PDFs and cubical persistence diagram for cases of Quintic Oscillator.}
        \label{fig:Quintic_analytical}
\end{figure}

Figure~\ref{fig:Quintic_analytical} shows the analytically generated PDFs for combinations of $h \in \{1, -1\}$ and $a \in \{1, 0, -1\}$, along with their super-level cubical persistence diagrams. Although the transition of $a \in [-1, 1]$ keeping $h=-1$ does not show a bifurcation, the shift from symmetry to asymmetry is an interesting phenomenon to observe and shows the system's ``preference" in terms of state. Note that for $a = 1$ and $a = -1$, the PDFs are topologically equivalent and discussing both of them does not add anything to the discussion on detecting bifurcations. Hence, the results for $a=1$ are being omitted for the estimated density case. For brevity, we will only discuss results for variation in $h$ with $a=0$ fixed and variation in $a$ with $h=1$ fixed.

\subsubsection{Varying $h \in [-1, 1]$ with $a=0$:}
\begin{figure}[!htbp]
\centering
\includegraphics[width=0.6\linewidth]{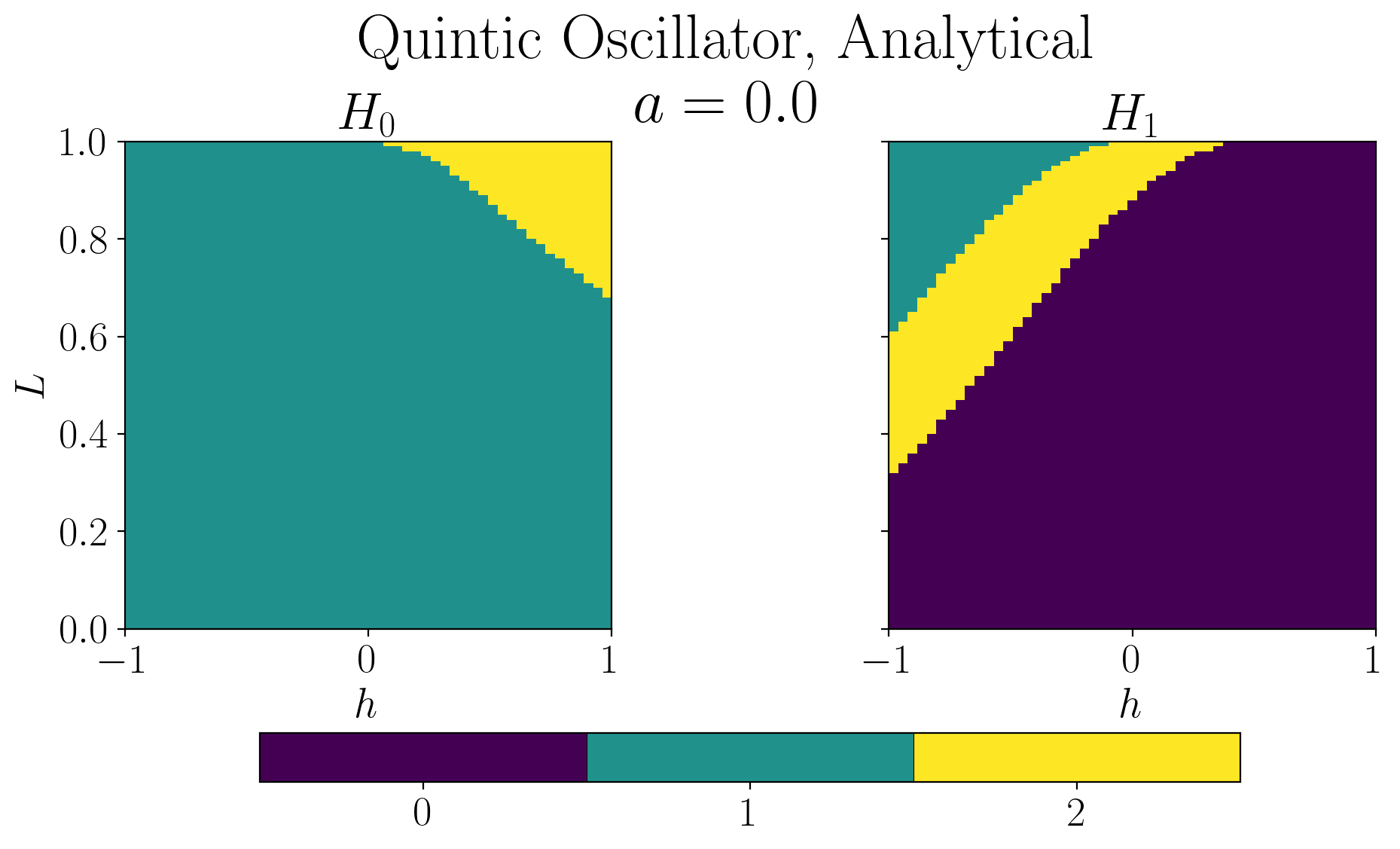}
\includegraphics[width = 0.6\textwidth]{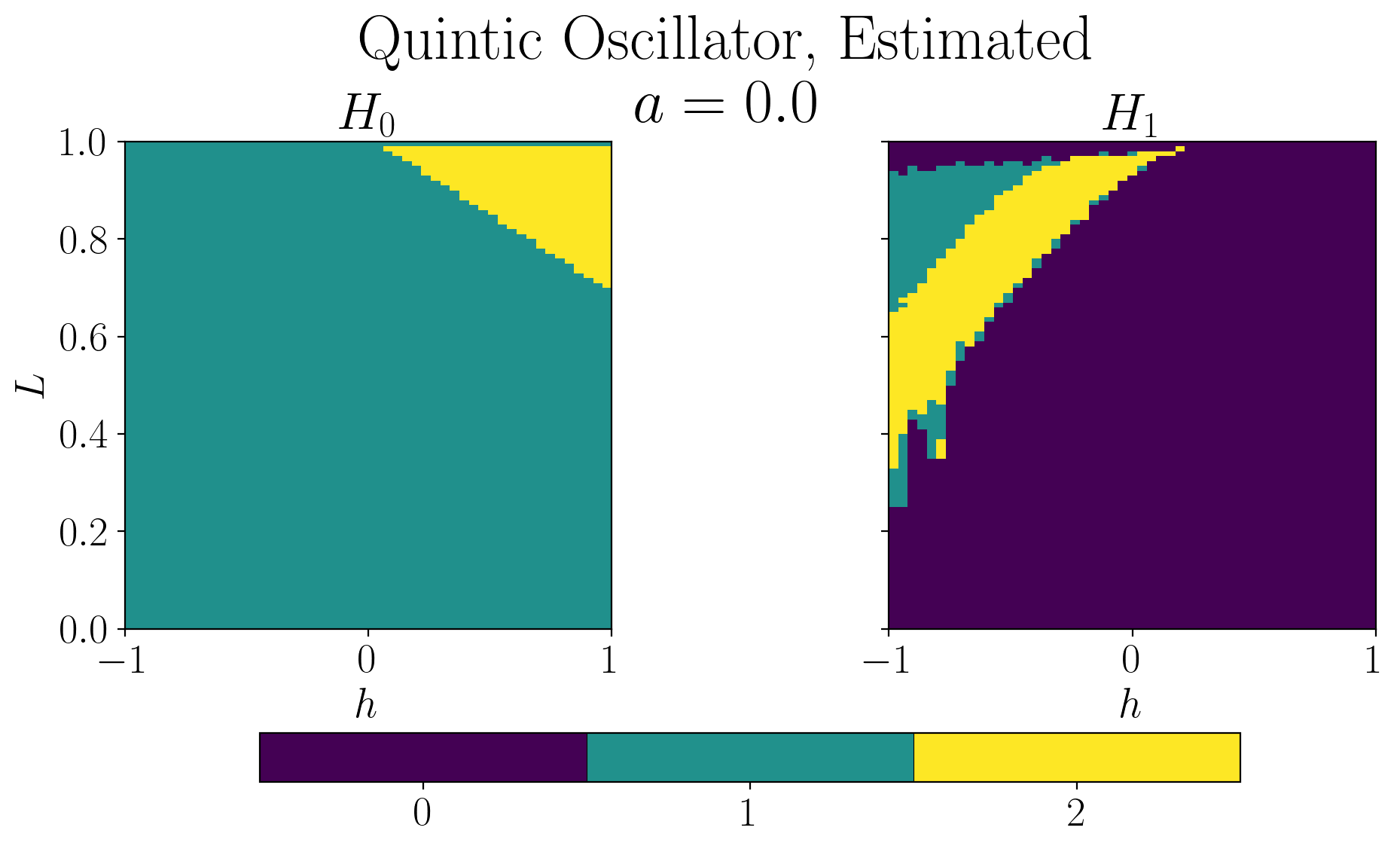}
\includegraphics[width=0.6\linewidth]{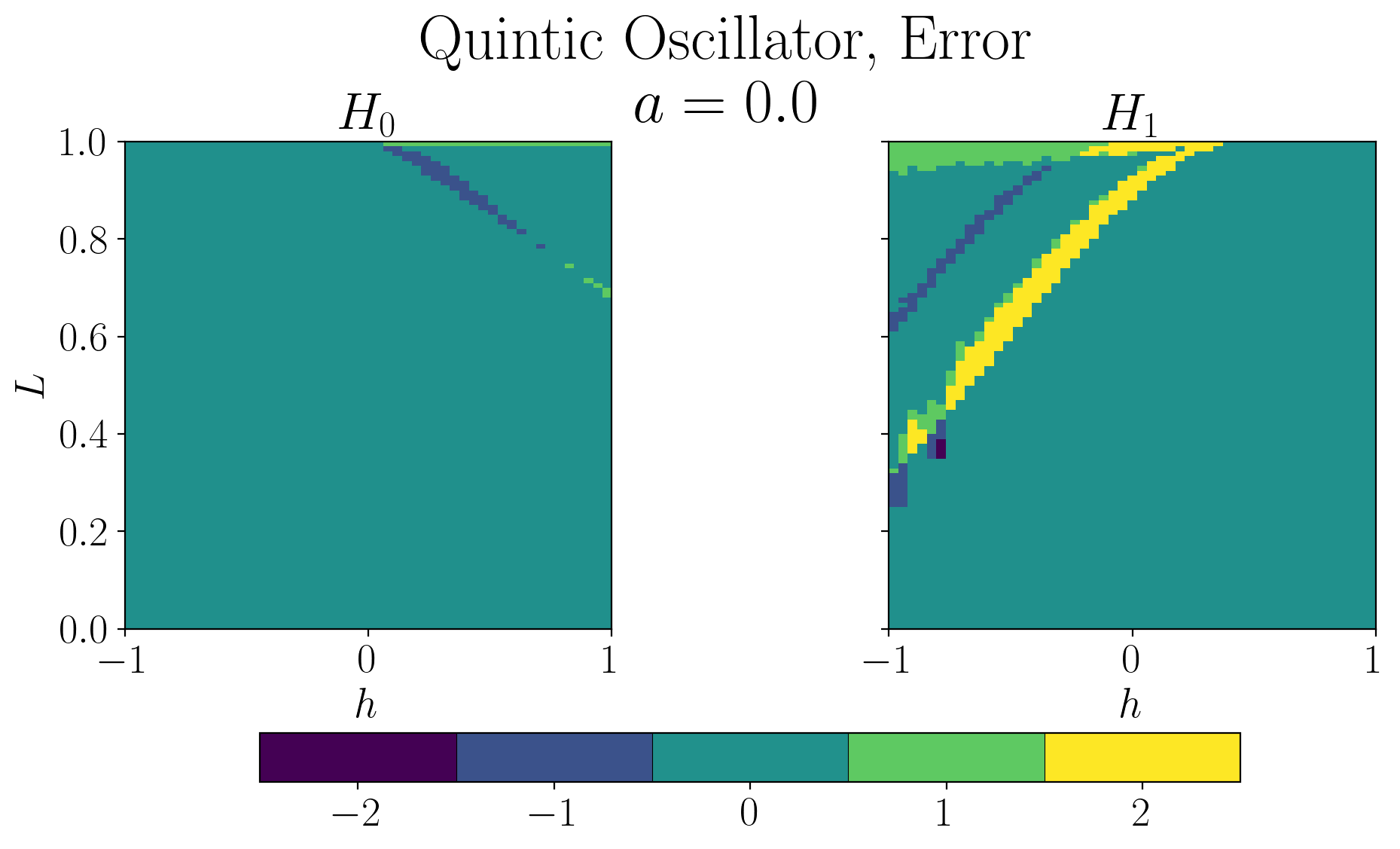}
\caption{Analytical and estimated homological bifurcation plots for Quintic Oscillator with $a=0$ fixed, with the errors.}
\label{fig:quintic_analytical_crocker_a0}
\end{figure}

Figure~\ref{fig:quintic_analytical_crocker_a0} shows the analytical and estimated homological bifurcation plots for Quintic Oscillator over the interval $h \in [-1, 1]$ keeping $a=0$, along with the errors. For $h>0$, the plots show the change in $H_0$ rank from $1$ to $2$ indicating the emergence of a second peak. For $h$ less than $\sim 0.5$, observe the change in $H_1$ rank indicating the emergence of inverted bi-stabilities (``holes") by the shift in rank from $0$ to $2$. This is accompanied by a further change in rank from $2$ to $1$, indicating the merger of the two voids into one as the level $L$ is raised. In general, the estimated plot has similar overall features as the analytical case, with most errors occurring around the topological boundaries, as expected.

\subsubsection{Varying $a \in [-1, 0]$ with $h=1$:}
Finally, we discuss the case for the interval $a \in [0, 1]$ while keeping the parameter $h=1$. Based on Fig.~\ref{fig:Quintic_analytical}, we expect the PDF to bifurcate from a monostable state with LCO (rank($H_0 = 2$), rank($H_1 = 1$)) to a bistable state (rank($H_0 = 2$), rank($H_1 = 0$)). 
\begin{figure}[!htbp]
\centering
\includegraphics[width=0.6\linewidth]{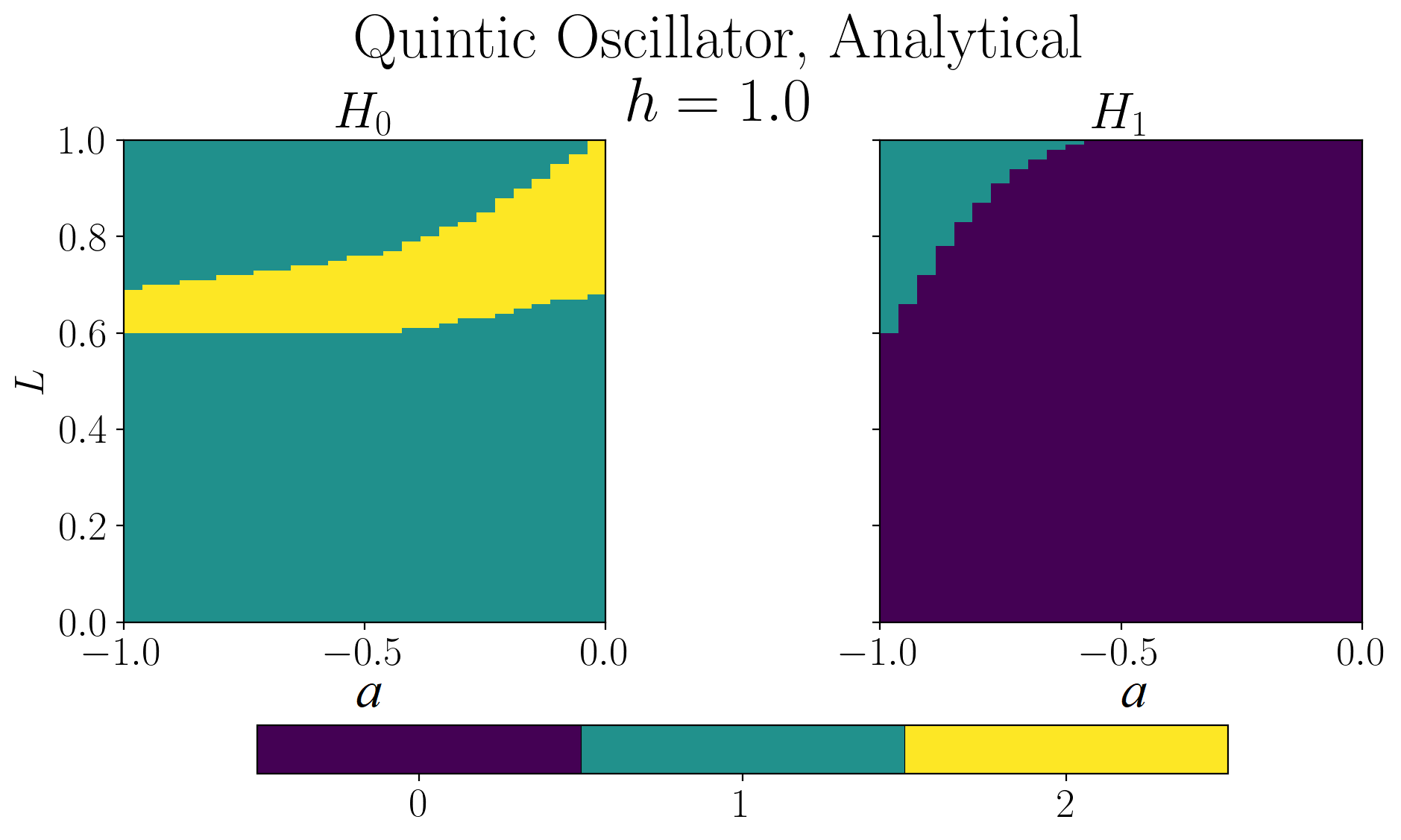}
\includegraphics[width = 0.6\textwidth]{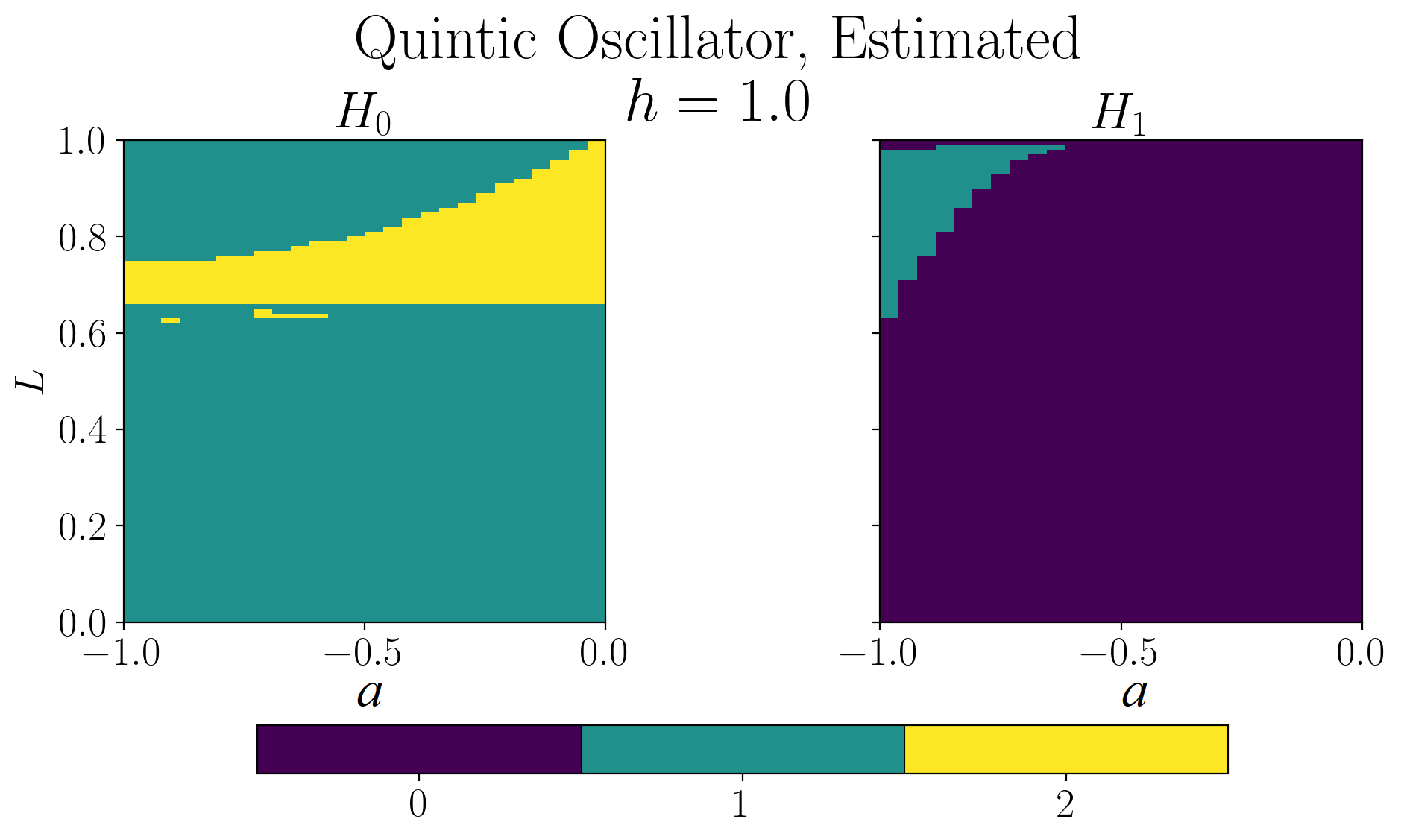}
\includegraphics[width=0.6\linewidth]{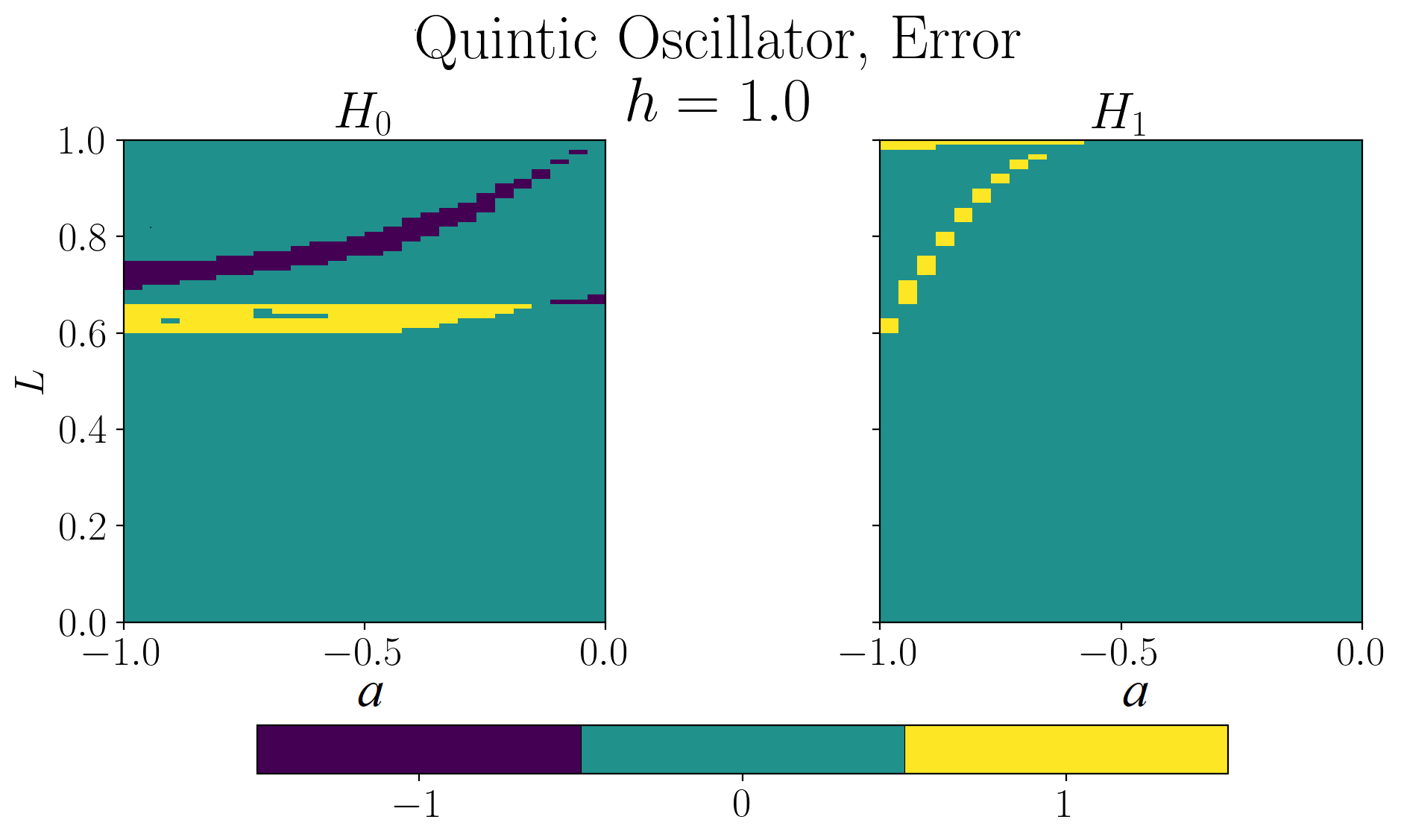}
\caption{Analytical and estimated homological bifurcation plots for Quintic Oscillator at $h=1$, with the errors.}
\label{fig:quintic_analytical_crocker_h1}
\end{figure}

Figure~\ref{fig:quintic_analytical_crocker_h1} shows the analytical and estimated homological bifurcation plots for Quintic oscillator with $h$ fixed at $1.0$, with the errors between them. For $a$ less than $\sim 0.55$, the $H_1$ plot displays the change in rank from $0$ to $1$ indicating the emergence of LCO. The $H_0$ plot displays an interesting thickness change in the rank $=2$ region as $a$ is varied. However, this thickness merely indicates an asymmetry in the heights of the two peaks (or one peak and one limit cycle) in the PDF. At $a=0$, the PDF's two peaks have the same maximum height, leading to a rank of $2$ throughout the level $L$. But as $a$ is varied, the peaks become asymmetric and do not share the same maximum height - leading to the change in homology rank from $2$ to $1$ at a level less than $1$. The estimated homology plots capture the same features as the analytical plots with errors occurring near the critical points.

\section{Conclusion}

Stochastic dynamical systems have wide applications in diverse fields. Accurately understanding, predicting and controlling the behaviour of these systems is of immense interest, and bifurcations (a qualitative change in the behaviour of the system) are a crucial aspect of such systems. P-bifurcations in stochastic systems are qualitative changes in the topology of the joint probability density function of the system response, due to which their analyses have often been restricted to qualitative methods. 

Despite the recent attempts to detect P-bifurcations, such as quantifying the critical points in the PDF, counting the number of peaks or computing the Shannon entropy over a range of bifurcation parameters, the methods and measures proposed are lacking in robustness and wide applicability. While some of these methods only work for one-dimensional systems, some are not general in their development and work for only specific systems. The remaining fail to capture richer dynamical behaviour or rely on the visual inspection of PDFs which makes it difficult to both estimate the stability bounds for the system and automate the procedure.

To address these limitations, our work provides a novel tool for detecting the qualitative P-bifurcation and understanding the state of the system using persistent homology, a tool from Topological Data Analysis (TDA). We illustrated that the change in the topology of the PDF (i.e. occurrence of a bifurcation) results in an abrupt change in the ranks of various homology classes, making it possible to quantify the onset of a bifurcation, and can be elaborately visualized using analytical homological bifurcation plots. Since analytical PDFs are rarely available for stochastic differential equations, we extended the concept to density estimates, by involving an algorithm for the topological consistency between analytical and estimated densities. While the algorithm for topological consistency introduces parameters such as $r$ and $\e$ which need to be carefully chosen, we demonstrated that the method is equipped to detect bifurcations using estimating homological bifurcation plots. The method was applied to three renowned oscillators: stochastic Duffing, Raleigh Vander-Pol (RVP) and Quintic, each exhibiting rich PDFs over their bifurcation intervals. We demonstrated that the method not only captured the Duffing oscillator's shift from monostable PDF to bistable PDF and the emergence of a stochastic limit cycle in RVP oscillator, but also detected the complex inverted bistable PDFs exhibited by the Quintic oscillator. Finally, it was demonstrated in each case that although discrepancies exist between the analytical and estimated plots, the errors are quite predictable and typically only occur along topological boundaries, as expected.

Since for this work, we had prior knowledge of the ground truth in each oscillator's case due to the availability of the analytical PDFs, choosing the optimal values for $r$ and $\e$ was comparatively easy to handle. However, in practice, this would be harder to accomplish and hence future work should be aimed at finding ways for automatic estimation of best choices for $r$ and $\e$. %
Other methods which reduce noise in persistence diagrams and provide topological consistency between analytical and estimated densities may be explored as well. Finally, an unspoken assumption throughout this work was the abundance of experimental data, which allowed for generation of meaningful density estimates. However, there may be instances where available samples from the system are singular (e.g. `Big Data') and may not result in statistically reliable density estimates. Therefore, future work should focus on developing quantitative bifurcation detection techniques for such elusive systems.

\section{Acknowledgements}
This material is based upon work supported by the Air Force Office of Scientific Research under
award number FA9550-22-1-0007.

\bibliography{SDE} 
\bibliographystyle{ieeetr}

\end{document}